\documentclass[11pt]{article}
\usepackage{xcolor}
\definecolor{blue}{rgb}{0,0,1}
\usepackage[pagebackref,bookmarks]{hyperref} 
\hypersetup{pagebackref=true,bookmarks=true,
       hyperfigures=true,
	unicode=true,
	colorlinks=true,
	citecolor=blue,
	linkcolor=blue,
	anchorcolor=red
}
\usepackage[textsize=footnotesize,color=blue!40, bordercolor=white]{todonotes}
\newcommand\todoi[1]{\todo[inline]{#1}}
\usepackage{etoolbox}
\makeatletter
\def\new@mathgroup{\alloc@8\mathgroup\mathchardef\@cclvi}
\patchcmd{\document@select@group}{\sixt@@n}{\@cclvi}{}{}
\patchcmd{\select@group}{\sixt@@n}{\@cclvi}{}{}
\makeatother

\usepackage{lmodern}
\usepackage{mathpazo}
\usepackage{fourier} %Best so far

\usepackage{amsmath}
\usepackage{amssymb}

\usepackage{graphicx,bbm,color,calc,ifthen,theorem,capt-of,tocloft}%,showkeys,showlabels}
\usepackage{amsmath}% needed before mathabx
\usepackage{amssymb}

\usepackage{fancyhdr}

\newcommand{\call}{ {\mathcal{L}}}

\newcommand{\phiabp}[2]{\ensuremath{\Phi_{\beta +1}^{(  \alpha )}}}

\newcommand{\nod}{{\noindent}}

\newcommand{\re}{\ensuremath{\mathbbm{R}}}

\newcommand{\rest}{{\upharpoonright}}
\newcommand{\emp}{{\varnothing}}
\newcommand{\pa}[1]{\ensuremath{\langle #1 \rangle}}

\newcommand{\la}{{\langle}}
\newcommand{\ra}{{\rangle}}
\newcommand{\ran}{\text{ran}}

\newcommand{\back}{{\backslash}}
\newcommand{\power}{P}

\newcommand{\ie}{{\itshape{i.e.}}{\hspace{0.25em}}}
\newcommand{\eg}{{\itshape{e.g.}}{\hspace{0.25em}}}
\newcommand{\etc}{{\itshape{etc.}}{\hspace{0.25em}}}
\newcommand{\via}{{\itshape{via }}{\hspace{0.25em}}}
\newcommand{\cf}{{\itshape{cf. }}{\hspace{0.25em}}}

\newcommand{\all}{{\forall}}
\newcommand{\Equi}{{\,\Longleftrightarrow\,}}
\newcommand{\equi}{{\,\longleftrightarrow\,}}
\newcommand{\ex}{{\exists}}
\newcommand{\Imp}{{\,\Rightarrow\,}}
\newcommand{\sset}{ { \subseteq } }

\newcommand{\imp}{{\,\longrightarrow\,}}

\newcommand{\rem}{{\noindent}{\bfseries{Remark: }}}
\newcommand{\fin}[1]{\ensuremath{[  ]^{< \omega}}}

\newcommand{\qed}{\hfill Q.E.D.}
\newcommand{\qedtwo}[1]{\hfill QED(#1)}
\newcommand{\pf}{{\noindent}{\textbf{Proof: }}}
\newcommand{\dfs}{\ensuremath{=_{\ensuremath{\operatorname{df}}}}}

\newcommand{\oddpagetext}[1]{\newcommand{\pageoddheader}{{\small }}}
\newcommand{\evenpagetext}[1]{\newcommand{\pageevenheader}{{\small }}}

\def\sset{\mbox{ $\subseteq$ }}
\catcode`\<=\active \def<{
\fontencoding{T1}\selectfont\symbol{60}\fontencoding{\encodingdefault}}
\catcode`\>=\active \def>{
\fontencoding{T1}\selectfont\symbol{62}\fontencoding{\encodingdefault}}
\catcode`\|=\active \def|{
\fontencoding{T1}\selectfont\symbol{124}\fontencoding{\encodingdefault}}

\newcommand{\nosymbol}{}
\newcommand{\nobracket}{}
\newcommand{\nocomma}{}

\newcommand{\tmem}[1]{{\em #1\/}}
\newcommand{\tmop}[1]{\ensuremath{\operatorname{#1}}}

\newcommand{\tmmathbf}[1]{${\boldmath{$ #1$}}$}
\newcommand{\tmtextit}[1]{{\it{#1}}}
\newtheorem{theorem}{Theorem}[section]

\newtheorem{corollary}[theorem]{Corollary}
\newtheorem{definition}[theorem]{Definition}

\newtheorem{lemma}[theorem]{Lemma}
\def\nod{\noindent}

\def\mb{\mbox{}}

\def\nod{\noindent}

\def\inn{\mbox{ {\large$\in$} }}

\def\ed{\end{document}}

\def\bth{\begin{theorem}}
\def\eth{\end{theorem}}

\def\blem{\begin{lemma}}
\def\elem{\end{lemma}}

\def\bdefn{\begin{definition}}
\def\edefn{\end{definition}}

\renewcommand{\qed}{\mbox{ } \hfill Q.E.D.}

\def\mathsf{}%%PDW Flashmode cannot cope with mathsf!
\def\cala{\mathcal{A}}
\def\cor{\color{red}}
\def\wq{\widetilde Q}
\def\wb{\overline{Q}}
\def\widebar{\overline}
\def\omicron{o}
\def\wmm{{\widetilde M}^{-}}

\def\trans{\tmop{Trans}}

\def\mid{\, | \, }
\def\siclub{$\tmmathbf{\Sigma}_{1}$-$\tmop{club}$}

\begin{document}

\title{Stably measurable cardinals}

\author{ P.D. Welch\\  School of Mathematics, \\University of Bristol,\\
  Bristol, BS8 1TW, England}

\date{December 27, 2018}

\maketitle

\begin{abstract}
  We define a weak iterability notion that is sufficient for a number of
  arguments concerning $\Sigma_{1}$-definability at uncountable regular
  cardinals. In particular we give its exact consistency strength firstly in terms of the second uniform indiscernible for bounded subsets of $\kappa$: $u_2(\kappa)$, and secondly to give the consistency strength of a property of L\"ucke's.
  
{\bf Theorem:} {\nod} The following are equiconsistent:\\ (i)  {\em There exists $  \kappa  $  which is stably measurable ;} (ii) {\em for some cardinal }$\kappa$, $u_2(\kappa)=\sigma(\kappa)$;\\
  (iii) {\em The {\boldmath$\Sigma_{1}$}-club property holds at  a cardinal }$\kappa$. \\  
\nod Here $\sigma(\kappa)$ is the height of the smallest 
$M \prec_{\Sigma_{1}} H ( \kappa^{+} )$ \ containing $\kappa+1$ and all of $H ( \kappa
  )$. % We give some  applications. 
  Let $\Phi(\kappa)$ 
 be the assertion: $ \forall X \sset \re\all r \in \re[
  X\mbox{ is }\Sigma_{1} ( \kappa, r )\mbox{-definable} \equi X\in \Sigma^{1}_{3}(r)].$
  
  %We obtain also:
{\bf Theorem: } 
%\begin{theorem}\label{2.1}
{\em  Assume $\kappa$ is stably measurable. Then $\Phi(\kappa)$.}\\
%the following are equivalent for
 % $X \sset \re$:
  %\quad\quad  (i) $X$ is $\Sigma_{1} ( \kappa )$-definable; \ \ (ii) $X$ %is
 % $\Sigma^{1}_{3}$ definable.}
And a form of converse:

{\bf Theorem: } {\em Suppose there is no sharp for an inner model with a strong cardinal. Then in the core model $K$ we have:\\
%{\em In $K^{strong}$ built below $0^{pistol}$, \\ 
$\mbox{``}\exists \kappa \Phi(\kappa) \mbox{'' is (set)-generically absolute} \equi$ 
%\quad \indent \quad 
$\mbox{There are arbitrarily large stably measurable cardinals.}$}

%In fact the existence of arbitrarily large stably measurable cardinals is consistently equivalent to ``$ \exists \kappa\Phi(\kappa)$'' being set-generically absolute.
When $u_2(\kappa) < \sigma(\kappa)$ we give some results on {\em inner model reflection.}
  
%  {\nod}The stable measurability property is much weaker than measurability, and in $K$ is equivalent to an assertion on the maximality of the length of the mouse order of  $H ( \kappa )$: that it be as long as the ordinal height of the least stable   set $M \prec_{\Sigma_{1}} H ( \kappa^{+} )$ \ containing all of $H ( \kappa
%  )$.
  
 % These properties also suffice to weaken the assumptions for some further results of L\"ucke, Schindler,
 % and Schlicht, including one on $F_{\kappa}$, the closed and unbounded filter on $\power(\kappa)$:
   
  % {\bf Theorem:} (ZFC) {\em If $\kappa$ is $\Sigma_{2}$-stably measurable, then $F_{\kappa}$ is not $\Pi_{1}^{H(\kappa^{+})}$-definable.}

\end{abstract}

\section{Introduction}

There are a number of properties in the literature that fall in the region of being weaker than measurability, but stronger than $0^{\#}$, and thus  inconsistent with the universe being that of the constructible sets.  Actual cardinals of this nature have been well known and are usually of ancient pedigree: Ramsey cardinals, Rowbottom cardinals, Erd\H os cardinals, and the like (\cf for example, \cite{Ka03}). Some concepts are naturally not going to prove the existence of such large cardinals, again for example, descriptive set theoretical properties which are about $V_{\omega+1}$ do not establish the existence of such large cardinals but rather may prove the consistency of large cardinal properties in an inner model. Weak generic absoluteness results, perhaps again only about $\re$, may require some property such as closure of sets under $\#$'s, or more, throughout the whole universe. 

An example of this is afforded by {\em admissible measurability} (defined below):\\
% (Def 3.5) \cite{W3} Theorem 4, Lemma 1

\nod {\bf Theorem}(\cite{W3} Theorem 4, Lemma 1) {\em Let $\Psi$ be the statement:
$$
\all D \sset \omega_{1}(D \mbox{ is universally Baire } \Equi \ex r \sset \omega ( D \in L[r] )).
$$
If $K$ is the core model then $\Psi^{K}$ is (set)-generically absolute if and only if there are arbitrarily large admissibly measurable cardinals in $K$.
}\\
%\end{theom}

This is a very weak property: weaker than an $\omega_{1}$-Erd\H os, but certainly stronger than {\em ``For any set $X\sset On,  X^{\sharp}$ exists'' }(thus indeed stronger than, say, two step $\Sigma^{1}_{3}$-generic absoluteness - see \cite{FMW}). Essentially it is often an assertion about the density of the mouse order in some, or alternatively arbitrarily large, $H(\kappa)$. This is also the guiding spirit behind the notions of stable measurability defined here.

In \cite{Luecke18} and  \cite{LSS18} the authors study, in essence, $\Sigma_{1}$-definable properties of a regular cardinal $\kappa$ in various forms: whether there is a $\Pi_{1}(\kappa)$ definition of the club filter on $\kappa$ for example, or whether $\Sigma_{1}(\kappa)$-definable subsets of $\kappa$ enjoy some kind of homogeneity property, such as that from \cite{Luecke18} defined below at \ref{sigmaclub}. The theorems of the abstract involve a strengthening of admissible measurability to stable measurability. This allows us an exact calibration of the strength of L\"ucke's {\boldmath$\Sigma_{1}$}-club property. It also allows minor improvements in the assumptions of certain theorems from $\cite{LSS18}$.

Stable measurability, whilst being ostensibly about $\Sigma_{1}$-definable subsets of $\kappa$, and whether an iterable measure can be put on the least stable set, is really something about the {\em bounded} subsets of $\kappa$. It says something about the strength of the mouse order in $H(\kappa)$ (the class of sets  hereditarily of cardinality less than $\kappa$), or relatedly, the size of the least uniform indiscernible above $\kappa$ for bounded subsets of $\kappa$. In the core model $K$, (at least below $0^{pistol}$) it is literally saying that the mouse order has length up to the least $\Sigma_{1}$ stable ordinal $\sigma(\kappa)$ as defined in this context. As the {\boldmath$\Sigma_{1}$}-club property turns out to be equivalent to stable measurability, it too, although phrased in terms of homogeneity properties of simply defined functions on $\kappa$, or subsets of $\kappa$, is in turn capable of being viewed as being actually about bounded subsets of $\kappa$.

Note: By {\em premouse} or {\em mouse} we mean that in the modern sense: see \cite{Z02}. By a {\em Dodd-Jensen mouse} (or $DJ$-mouse) we mean that of \cite{D}. We do not need many details of the latter: simply that they are similar to the levels of $L[\mu]$ where the levels are defined as in simple relativised constructibility from a predicate $\mu$. A $DJ$-mouse then is a structure of the form $\la J_{\alpha}^{U},\in, U\ra\models \mbox{ ``  $U$ is a normal measure on $\kappa $''}$ with wellfounded iterated ultrapowers. Another required feature of a $DJ$-mouse $M$ is that there is always a new subset of the measurable cardinal $\kappa$ definable over $M$. Consequently there is always also a definable onto map $f:\kappa \imp J_{\alpha}^{U}$. The Dodd-Jensen core model $K^{DJ}$ can be thought of as an $L[E]$ hierarchy whose initial segments are all sound mice in the usual fashion, or alternatively as simply the union of the older $DJ$-mice. These universes are the same. Whenever the $K^{DJ}$ model is mentioned, for fixity we shall assume the former, modern, now standard, presentation.

In $K^{DJ}$ there is a natural method of comparison of $DJ$-mice in $H(\kappa)$: iterate them all $\kappa$ times, and the union thus obtained is the ``$Q$-structure at $\kappa$'', and is of the form $Q = Q(\kappa)=(J^{F_\kappa}_\theta, \in, F_\kappa)$ for some ordinal $\theta(\kappa)$ where $F_\kappa$ is the cub filter on $\kappa$, but which is an amenable iterable measure   on $Q$. This is a useful structure to work with even if it does not conform to the modern notion of mouse. With sufficiently many $\sharp$'s in $H(\kappa)$, $\theta(\kappa)$ can be (but is not always) the second `uniform indiscernible' for bounded subsets of $\kappa$.  But if it is then the critical points of the iterates of $Q$ enumerate precisely these uniform indiscernibles.

All of this is in $K$ the core model. However here in this paper we also step out of $K$ and look at generalizations $\wb(\kappa)$ (Def.\ref{1.14}) and similar characterisations that now generate the uniform indiscernibles in $V$. Roughly speaking the greater the ordinal height of $\wb(\kappa)$ (corresponding to the earlier ordinal height $\theta(\kappa)$ of $Q(\kappa)$ in $K$) the `stronger' the iterability properties instantiated in $H(\kappa)$. 

If we approach from the other direction and ask if any subsets of $\kappa$ (rather than bounded subsets of $\kappa$) can be put in sufficiently closed iterable structures $(M,\in, U)$ (think of putting any subset of $\kappa$ in a transitive $\kappa$-sized models $M= \mb ^{<\kappa}M$ with a wellfounded ultrapower map $j:M\rightarrow N$ to get weak compactness) then we get a notion of {\em iterable cardinal}. This is of course weaker than measurability, but it is also weaker than Ramseyness (\cite{Sharpe-Welch2012} Lemma 5.2) which requires (as Mitchell \cite{Mi79}, Jensen \cite{DJK} showed) not just that $(M,\in, U)$ be iterable but that additionally $U$ be $\omega$-closed.

Several of the theorems of \cite{Luecke18} , \cite{LSS18} use as an iterability assumption that $\kappa$ be an iterable cardinal.  We observe here that instead one needs only something weaker: that a $\Sigma_1$-substructure $N$ of $H(\kappa^+)$ be itself placed in such an iterable $(M,\in, U)$. This is the notion of being $(\Sigma_1)${\em-stably measurable}.  That this is not just some minor improvement resides in the fact that some of the properties turn out to be equiconsistent to stable measurability, or even equivalent in a canonical inner model such as $K^{DJ}$.

\setcounter{section}{2}

\setcounter{theorem}{5}
\begin{theorem}%\ref{ordinalsinK}
$(V=K^{DJ})$ 
$$\sigma ( \kappa ) =u_{2} ( \kappa )   \equi \kappa\mbox{ has the
  {\boldmath{${\Sigma}_{1}$}}-}\tmop{club}\mbox{ property }\equi \kappa \mbox{ is stably
  measurable. }$$
\end{theorem}

\setcounter{theorem}{0}
\setcounter{section}{3}

Our theorem in the analogous form to that which began this introduction is spread over the following two statements.
 We have:

\begin{theorem}%\ref{2.1}
Let $\Phi(\kappa)$ be the following sentence:
$$ \Phi(\kappa) : \forall X \sset \re\all r \in \re[
  X\mbox{ is }\Sigma_{1} ( \kappa, r )\mbox{-definable} \equi X\in \Sigma^{1}_{3}(r)].$$
  Assume $\kappa$ is stably measurable. Then $\Phi(\kappa)$ holds.
  
  %the following are equivalent for   $X \sset \re$:
  
%  (i) $X$ is $\Sigma_{1} ( \kappa )$-definable; \ \ (ii) $X$ is
%  $\Sigma^{1}_{3}$ definable.
\end{theorem}

In one sense we have an equivalence:
\setcounter{theorem}{2}
\begin{theorem}%\ref{2.3} 
Assume $V=K^{DJ}$.
% and that $H(\kappa)$ is closed under sharps.  
 $$\kappa\mbox{ is stably measurable }\equi  \Phi(\kappa) \mbox{ is preserved by small forcings of size } < \kappa.$$

\end{theorem}

\setcounter{theorem}{5}

\begin{corollary} Assume $V=K^{DJ}$ (or $V=K^{strong}$). % and is closed under sharps. 
Then $\exists \kappa \Phi(\kappa)$ is (set)-generically absolute if and only if there are arbitrarily large stably measurable cardinals.
\end{corollary}

\setcounter{theorem}{0}
\setcounter{section}{1}

Our theme in essence is to tease out the implications between the notions of stable measurability,  good  {\boldmath$\Sigma_1$}$(\kappa)$-wellorders, and the length of the mouse order when working in $L[E]$ models, or, when in $V$, the height of the $\wb(\kappa)$-structure  which contains all the $\kappa$'th  iterates of coarse `mouse-like' objects in $H(\kappa)$. 

In the final section we make some comments on {\em inner model reflection} by identifying the least $L[E]$ models which reflect $\Pi_{n}$ sentences into their inner models. Such a model is then not `pinned down' by such a sentence (with ordinal parameters allowed). This phenomenon occurs before stable measurability, and can be seen to happen when $u_{2}(\kappa)<\sigma(\kappa)$, but the mouse order is sufficiently long to be beyond `admissible measurability'.

%\begin{definition}.   \end{definition}
\subsection{Stable Measurability}
\begin{definition}\label{kappamodel} We say that $N$ is a $\kappa$-{\em model} if: $\trans(N)$, $\kappa \in N$ and $\mb ^{<\kappa}N\sset N$.

\end{definition}

\begin{definition}
  \label{sm} Let $\omega < \kappa \in \tmop{Reg}$. Then $\kappa$ is
  $\Sigma_{n}$-{\tmem{stably measurable}} if, for some transitive $M
  \prec_{\Sigma_{n}} H ( \kappa^{+} )$ with $M \supseteq H ( \kappa ) \cup \{
  \kappa \}$, \ there is a $\kappa$-model $N\supseteq M$ and a filter $F$ with $( N, \in ,F ) \models$``$F$ is a
  normal measure on $\power ( \kappa )$'' so that $( N, \in ,F )$ is amenable, and it is iterable,
  that is, has wellfounded ultrapowers by the measure $F$ and its images. We say that 
 ($M$ and) $( N, \in ,F )$ {\em ``witnesses $\Sigma_n$-stable measurability.''}
\end{definition}

%\todoi{Amenability should go here?}
The above is by way of analogy with the notion of {\tmem{admissibly
measurable}} which was coined in {\cite{W3}}. This required only that $M$ be
the least transitive admissible set containing $H ( \kappa ) \cup \{ \kappa
\}$ and again with an appropriate filter $F$ with wellfounded ultrapowers. In
the above if $n=1$ we just refer to stable measurability.

\begin{definition} We say that $\prec$ is a {\em good $\Sigma_{1}(p)$-wellorder} of $\power(\kappa)$ if $\prec$ as a binary relation has a $\Sigma_{1}(p)^{H ( \kappa^{+} )}$ definition (in some parameter $p\in H(\kappa^{+})$, and so that the set of all initial segments $\{z \mid \ex x\in \power(\kappa) \wedge z = \{y\mid y \prec x\}\}$ is a $\Sigma_{1}(p)^{H ( \kappa^{+} )}$ set.
\end{definition}

Note: (i) if there is a good $\Sigma_{1}^{H ( \kappa^{+} )}(\{\kappa, p\})$ wellorder of
$\power ( \kappa )$, (for some $p\in H(\kappa)$) we can define $\Sigma_{1}$-Skolem functions in the usual
manner and more readily define such an $M$. In some $L [ E ]$ models this will
be the case, and we shall use below the example of the Dodd-Jensen core model
$K=K^{\tmop{DJ}}$.

(ii) For $\Sigma_1$-stability ($n=1$) we shall show that we can take $N$ as an $M$ which is itself a $\Sigma_1$ elementary substructure.
If $\pa{N, \in ,F}$ witnesses stable measurability at $\kappa$, we
should just emphasise that without additional requirements, we cannot assume
that it is an iterable premouse of any form of the usual definition(s) of
premouse.

(iii) If $\kappa$ is $\Sigma_1$-stably measurable then it is easily seen to be a Mahlo cardinal. (If there is a $C\sset \kappa$  a cub set of singular cardinals, then there is such in $M\prec_{\Sigma_{1}}H(\kappa^{+})$. Now as $M$ is in  some iterable $N$ if $j:N\imp N'$ is the first ultrapower of $N$ by the $N$-normal measure, then $\kappa \in j(C)$ is singular in $N'$ which leads to a contradiction.)

(iv) Just using the increased elementarity available it is easy to see that for any $n\geq 2$ that $\Sigma_{n}$-stable measurability is equivalent to iterability. Hence we shall mostly be interested in
 $\Sigma_{1}$-stable measurability (and drop the ``$\Sigma_{1}$'').

\begin{definition}\label{sigma} We set
  $\sigma = \sigma ( \kappa ) = \tmop{On} \cap M$ to be the least ordinal which is the height of a transitive $M$ with $M \prec_{\Sigma_{1}} H (
  \kappa^{+} )$ and $M \supseteq H ( \kappa ) \cup \{ \kappa \}$.
\end{definition}

We shall remark below that our definition of stable measurability will ensure that there is such an $M$ as a least
$\Sigma_{1}$-substructure of $( H ( \kappa^{+} ) , \in )$ containing $H (
\kappa ) \cup \{ \kappa \} ,$ even in the absence of some canonical wellorder,
or canonically chosen skolem functions, for \ $H ( \kappa^{+} )$.

\begin{definition}
  Let $M_{0}=M_{0}(\kappa) \dfs \left\{ A \mid  A \sset \kappa \wedge \{ A \} \right.$ is a
  $\Sigma_{1} ( \kappa ,p )$-singleton set for some $\nobracket p \in H (
  \kappa ) \}$.
\end{definition}

In the above we could have written $\{ A \}$ \ is to be a $\Sigma_{1}^{H (
\kappa^{+} )} ( \kappa ,p )$-singleton, by Levy-absoluteness.

\begin{definition}
(i)  For $A \sset \kappa$ let $\sigma_{A} \dfs$ the least $\sigma > \kappa$ such
  that $L_{\sigma} [ A ] \prec_{\Sigma_{1}} H ( \kappa^{+} )^{L[A]}$. \\ (ii) 
  $\widetilde M=\widetilde M(\kappa) \dfs \bigcup_{A \in M_{0}} L_{\sigma_{A}} [ A ]$. \\
  (iii) $   \wmm= \wmm(\kappa)  \dfs \bigcup \{ L_{\sigma_{a}} [ a ]\mid {a {\sset} \gamma <\kappa,\,
  a^{\sharp} \mbox{ exists}}\}.$
  %then (i) $( M, \in )  \prec_{\Sigma_{1}} ( H ( \kappa^{+} ) , \in )$; (ii) $( M, \in )$ is admissible.
\end{definition}

The last definitions might seem peculiar at first glance, but they are suitable for analysing certain sets when we do not assume a 
good {\boldmath$\Sigma_1$}$(\kappa)$-wellorder of $\power(\kappa)$.    $\widetilde M$ can be thought of as  an approximation to a $\Sigma_{1}$-substructure of $H(\kappa^{+})$. Add a 
good {\boldmath$\Sigma_1$}$(\kappa)$-wellorder and it will be (see Lemma \ref{M} below). Moreover stable measurability of $\kappa$ will imply 
(Lemma \ref{key}) that $\widetilde M = \wmm$. It is this last equality that prompts the idea that $\Sigma_{1}$-stability of $\widetilde M$ is really about the bounded subsets of $\kappa$.

\begin{lemma}\label{allaresingletons}  Every $x \in \widetilde M$ is coded by some $B \in M_{0}$.
\end{lemma}
{\pf} Fix an $x \in \widetilde M$; there is thus some $A \in M_{0}$,  $\alpha < \sigma_{A}$, with $x
\in L_{\alpha} [ A ]$.
 Standard reasoning shows that there are arbitrarily large
$\beta < \sigma_{A}$ with $J^{A}_{\beta +1} \models$``$\kappa$ is the largest
cardinal'' and so that there is a $\Sigma^{J^{A}_{\beta +1}}_{1}( A,
\kappa )$ definable function $f: \kappa \imp J^{A}_{\beta +1}$.  We may  further assume
that $T,$ the $\Sigma_{1}$-$\tmop{Th} ( J^{A}_{\beta +1} , \in ,A )$ coded as
a subset of $\kappa$ is in fact a $\Sigma_{1} ( \kappa ,A,q )$ singleton, for
some $q \in H ( \kappa )$, and hence a $\Sigma_{1} \left( \kappa ,  \la p,q
\ra \right)$-singleton where $\{ A \} \in \Sigma_{1} ( \kappa ,p )$. (This is
because we can take $T$ as the unique $\Sigma_{1}$-Theory of a level in the $L
[ A ]$ hierarchy where some $\Sigma_{1}$ sentence $\psi ( q )$ about some $q
\in L_{\kappa} [ A ]$ first becomes true.) But then from the theory $T$ we obtain $f$ and then may define
$\pa{\tmop{TC} ( \{ x \} ) , \in} \cong B_{0} \dfs \left\{ \pa{\xi_{0} , \xi_{1}}
\mid f ( \xi_{0} ) \in f ( \xi_{1} ) \in f ( \zeta ) \right\}$ for some $\zeta
< \kappa$, if $\{ x \} \in J^{A}_{\beta +1}$. Coding $B_{0}$ by G\"odel pairing as subset of $\kappa$, $B$,  we have $\{ B \} \in \Sigma_{1}
\left( \kappa , \pa{p,q, \zeta} \right)$ and so $B \in M_{0}$ as required. 
\ \ {\qed}

\begin{lemma}
  \label{M} Suppose there is a good $\Sigma_{1}^{H(\kappa^{+})}(\kappa,p)$ wellorder of $\power(\kappa)$ for some $p\in H(\kappa)$. Then $\widetilde M \prec_{\Sigma_{1}}H(\kappa^{+})$.
\end{lemma}
\pf
Using the good wellorder we have $\Sigma_{1}$-skolem functions for $\la H(\kappa^{+}),\in \ra$ which are themselves ${\Sigma_{1}}^{H(\kappa^{+})}$. Suppose that we have for each  $\Sigma_{1} \,\exists v_{0}\varphi(v_{0},v_{1})$ a skolem function $f_{\varphi}$ so that for all $A\subseteq \kappa$ if there is $u$ so that $\varphi(u,A)$ then $H(\kappa^{+}) \models \varphi(f_{\varphi}(A),A)$ holds. Suppose that $\exists v_{0}\varphi(v_{0},A)$ holds with $A\in M_{0}$. Then we may assume that the witness $u$ is itself a subset of $\kappa$ which is a $\Sigma_{1}(\kappa, A)$ singleton. This is because every set in $H(\kappa^{+})$ has cardinality there less than or equal to $\kappa$; given the good wellorder, we thus have for every $u$ there is a least, in the sense of the wellorder, subset of $\kappa$, $U$ say, that codes a $u$ that witnesses $\varphi(u,A)$. Then $\{U\}$ is a $\Sigma_{1}(\kappa, A, p)$-singleton, and so  $U \in M_{0}\sset \widetilde M$. Putting this together we have that $(\exists v_{0}\varphi(v_{0},A))^{\widetilde M}$. \qed\\
%\todoi{Check the pararemter p used in the last proof.}

\subsection{On $ \widetilde Q$}

\begin{definition}\label{1.14}
Let $\wq(\kappa)$ denote:
$$ \bigcup \{ N_{\kappa}\mid N_{\kappa} \mbox{ is the } \kappa \mbox{'th iterate of some amenable iterable } \pa{N,\in, U}\in H(\kappa)\} 
$$
\end{definition}

Under the hypothesis of the next lemma we shall have that $\wq$ is rud. closed.

\begin{lemma}\label{Lemma1.10} Suppose all bounded subsets of $\kappa$ have sharps. Then $\wmm(\kappa) = \widetilde Q(\kappa) $. Additionally any $X\in \power(\kappa)^{\wq}$ either contains or is disjoint from a set cub in $\kappa$.
\end{lemma}
\pf $(\sset)$: If $x \in \wq$ then for some $a = \la N,\in U\ra \in H(\kappa)$, $x \in N_\kappa$. But $N_\kappa\in L_{\sigma_a}$. So $x \in \wmm$.\\
$(\supseteq)$: Let  $x\in L_{\sigma_a}[a]\cap \power(\kappa)$ some $a \in H(\kappa)$. As $a^\sharp$ exists, let $N^a$ be the $a^\sharp$ mouse. Then $L_{\kappa'}[a]\sset (N^a
)_\kappa$ where $\kappa'= (\kappa^+
)^{L[a]}$. As $\sigma_a < \kappa'$, $x\in (N^a
)_\kappa \sset \wq$.  This shows that any such $x$ will be disjoint from, or contain a tail of the cub set of the sequence of iteration points of $N^a$.
\qed
%\todoi{Then in $K^{DJ}$ this should be the same as $Q$; if $H(\kappa)$ closed under sharps, then in V this should be the same as $Q$??; if $\kappa$ s.m. then $\widebar Q$ should be $\widetilde M$. }
%\begin{corollary} Suppose all bounded subsets of $\kappa$ have sharps. 

%\end{corollary}

\begin{lemma} Suppose all bounded subsets of $\kappa$ have sharps. Then (i) $\widetilde Q$ is rudimentary closed; (ii) $\pa {\widetilde Q,F_{\kappa}}$ is amenable and iterable, with $F_{\kappa}\cap \widetilde Q$ a $\widetilde Q$-normal ultrafilter.

\end{lemma}
\pf (i) As the rudimentary functions have as a generating set a finite set of binary functions (\cite{Je72}), it suffices by the last lemma, since each $L_{\sigma_{a}}[a]$ is rud. closed (it is an admissible set), to show that if $X,Y\in \widetilde Q$, that there is $c$ a bounded subset of $\kappa$ with $X,Y\in L_{\sigma_{c}}[c]$. By our supposition any  $a\in H(\kappa)$ is a member of the least $a$-mouse generating $a^{\sharp}$, $N^{a}$, and moreover 
$L_{\sigma_{a}}[a]\in N^{a}_{\kappa}$, the $\kappa$'th iterate of $N^{a}$. But then it is trivial that if $\{X\}\in N^{a}$ and $\{Y\}\in N^{b}$ then $\{X,Y\}\in N^{a\oplus b}\in \widetilde Q$ as $L_{\sigma_{a}}[a] \cup L_{\sigma_{b}}[b]\sset L_{\sigma_{a\oplus b}}[a\oplus b]$.

For (ii): %Any $A\in \power(\kappa)\cap \tilde Q$ is in some $L_{\sigma_{a}}[a]$ for an $a$ bounded in $\kappa$. As indicated any such $A\in N^{a}_{\kappa}$ where $N^{a}_{\kappa}$ is a mouse on $\kappa$ with its measure generated by the final segment filter of its critical points below $\kappa$. Thus every such $A$ is either disjoint from, or contains a c.u.b. set in $\kappa$. 
That $F_\kappa$ measures $\power(\kappa)\cap \widetilde Q$ is the last corollary.
%This shows $\pa {\widetilde Q,F_{\kappa}}\models $``$F_{\kappa}$ is an ultrafilter''. 
For amenability just note that any $\pa {Z_{\nu}\mid \nu <\kappa}\in L_{\sigma_{a}}[a]$ is again in  $\pa {N^{a}_{\kappa}, F_{\kappa}\cap N^{a}_{\kappa}}$. But the latter structure is amenable, (this is true of any $a$-mouse) and so $\{\nu \mid Z_{\nu}\in F_{\kappa}\}\in N^{a}_{\kappa}\in \widetilde Q$. Normality of 
$F_{\kappa}\cap \widetilde Q$ in $\widetilde Q$
is similar, and iterability follows from the countable closure of $F_{\kappa}$. \mb \hfill
\qed \\

For notation we set $I^{c}$, the  closed and unbounded class of Silver
indiscernibles for $L [ c ]$, to be enumerated as $\pa{\iota_{\alpha}^{c} \mid 
\alpha \in \tmop{On}}$ for $c$ a set of ordinals. 
\begin{definition}
  Suppose for every bounded subset $b $ of $\kappa \nocomma \nocomma  $,
  $b^{\#}  $exists. Then set
   $$u_{2} ( \kappa ) = \sup \{
  \iota^{b}_{\kappa +1} \, |\, b  \mbox{ a bounded subset of } \kappa
  \} .$$
  More generally: \,
   $$\pa{ u_{\iota} ( \kappa )\mid 0<\iota \in On}$$
   % = \sup \{
  %\iota^{b}_{\kappa +1} | b \in H ( \kappa ) \cap \power ( \kappa )\} 
  enumerates in increasing order $\bigcap \{ I^{b}\mid b \mbox{ a bounded subset of } \kappa\}$  
\end{definition}

%\todoi{Here we don't need the good $\Sigma_{1}$-WO}

Then this is by way of analogy for the second uniform indiscernible for the
reals, but now for bounded subsets of $\kappa$. By the same arguments as for
reals, $u_{2} ( \kappa )$ is also $\sup \{ \kappa^{+L [ b ]} \mid b \in
H ( \kappa ) \cap \power ( \kappa ) \} $. Indeed, as is well known, for any successor $\iota +1$: $$u_{\iota +1}= 
\sup \{ u_{\iota}^{+L [ b ]} \mid b  \mbox{ a bounded subset of } \kappa
  \} = 
\sup \{ \iota^{b}_{u_{\iota}+1} \mid b  \mbox{ a bounded subset of } \kappa
  \} .$$

It is an exercise in the use of sharps to add to this  that 
$u_{2}(\kappa)= 
\sup \{ \sigma_b \mid b  \mbox{ a bounded subset of } \kappa
  \}$.
The size of $u_{2}(\kappa)$
with reference to $\kappa$, gives, roughly speaking, the length of the mouse
order on $H ( \kappa )$. Indeed in $L [ E ]$ models (at least below a strong
cardinal) this can be made precise. Thus the next lemma interpreted in for
example, the Dodd-Jensen core model $K^{\tmop{DJ}}$, is declaring the length
of the mouse order restricted to $H ( \kappa )$ there, as somewhat long. In fact it will turn out to be maximal for this model.

\begin{lemma} Suppose that $H(\kappa)$ is closed under sharps. Then  the critical points of the iterated ultrapowers of $\pa {\widetilde Q,F_{\kappa}}$ are the uniform indiscernibles $\pa {u_{\iota}(\kappa)\mid 0 < \iota\in On}.$ Moreover if $\pa  {\widetilde Q_{\alpha},F^{\alpha}}_{\alpha\in On}   $ is the iteration of $\pa  {\widetilde Q_{1},F^{1}} = \pa {\widetilde Q,F_{\kappa}}$, with iteration maps $j_{\alpha,\beta}\, (1\leq \alpha < \beta \in On)$, and critical points $\lambda_{\alpha}(1\leq  \alpha \in On)$ then

\begin{eqnarray} u_{\alpha}(\kappa) &=& \lambda_{\alpha}\\
 u_{\alpha +1}(\kappa) & = &\widetilde Q_{\alpha} \cap On 
 \end{eqnarray}

\end{lemma}

\pf  First we note that as $\pa {\widetilde Q,F_{\kappa}} = \pa  {\widetilde Q_{1},F^{1}}$ is a rudimentary closed structure, we can prove a Los Theorem for its ultrapowers and the usual result for such a structure that it is a $\Sigma_{0}$ preserving embedding which is {\em cofinal} (that is if $k:\pa {\widetilde Q,F_{\kappa}}\imp Ult (\widetilde Q,F_\kappa)$, and if $\pi :Ult (\widetilde Q,F_\kappa)\imp (\widetilde Q_{2},F^{2})$  is the transitive collapse map, then taking $j= j_{1,2}=\pi \circ k$ we have that 
$\forall x\in \widetilde Q_2\exists y \in \widetilde Q ( x\sset j(y))$). Thus $j$ is in fact $\Sigma_{1}$-preserving. Note that by the amenability of $\pa {\widetilde Q,F_{\kappa}}$, $\power(\kappa)\cap \widetilde Q = \power(\kappa)\cap \widetilde Q_{2}$.  Suppose now  $[f] < [c_{\kappa}]$ in 
$Ult (\widetilde Q,F_\kappa)$. Thus
$f\in \widetilde Q, \, f:\kappa \imp On \cap \widetilde Q$ and by normality, with $\{\xi \mid f(\xi)< \kappa \}\in F_{\kappa}$.  Thus for $a \in H(\kappa)\cap \power(\kappa)$ we shall have $f\in L_{\sigma_{a}}[a]$.  By Silver indiscernibility $f(\xi) = h^{L[a]}(i,a, \vec \gamma, \xi , \vec \gamma')$ for some $\vec \gamma, \vec \gamma' \in [I^{a}]^{{<\omega}}$ with $\tmop{max}(\vec \gamma) \leq \xi \tmop{min}(\vec \gamma'))$ and $h^{L[a]}$ a canonical $\Sigma_{1}$-skolem function for $(L[a],\in, a)$. But going to $a^{\sharp}$ we shall have 
$f(\xi) = h^{L[a^{\sharp}]}(i',a, \vec \gamma, \xi )$ for some $i'$. In particular $f(\xi)<\gamma' =g(\xi)\dfs \tmop{min} I^{a^{\sharp}}\backslash \tmop{max}(\vec \gamma, \xi) +1$. Let $\gamma_{0}' \dfs 
\tmop{min} I^{a^{\sharp}}\backslash (\kappa
 +1)$. Then $\gamma_{0}' < On \cap \wb _{1}$. But then we have that $[f] < [g]$ and $j(f)(\kappa)<j(g)(\kappa)< \gamma_{0}' < On \cap \wb _{1}$. This shows that $j(\kappa) \leq On \cap \wb _{1}$.  But clearly as well $j(\kappa) \geq On \cap \wb _{1}$.

Thus (recalling that $\lambda_{1} = \kappa$ and $ \widetilde Q_{1}=\wq $):
$$u_2(\kappa) = \tmop{sup}\{\iota^{a}_{\lambda_{1}+1}\mid a \in H(\kappa)\cap \power(\kappa)\} = 
\tmop{sup}\{\sigma_{a}\mid a \in H(\kappa)\cap \power(\kappa)\}
= On \cap \widetilde Q_{1} .$$
But we have just seen that $j_{1,2}(\lambda_{1})=\lambda_{2}= On \cap \widetilde Q_{1}.$ This establishes (1) for $\alpha = 2$, and (2) for $\alpha = 1$, and the reader can deduce the cases for larger $\alpha$ from this.\qed\\

This then gives a simple expression for the uniform indiscernibles of the bounded subsets of $\kappa$: they are the iteration points of $ \pa {\widetilde Q,F_{\kappa}} $ as well as (their successor) elements being  the ordinal height of the ultrapowers.  (The reader will recall that under $AD$, in $L(\re)$ we have that for reals, $u_{2}=\aleph_{2}$ and the ultrapower of $\pa{u_{1},<}/ F_{\omega_{1}}$ is $u_{2}$.)
The following is well known for reals  but follows immediately from the above: 
\begin{corollary}$\tmop{cf}(u_{\alpha+1}(\kappa)) =\tmop{cf}(u_{2}(\kappa)) $.
\end{corollary}
\pf $j_{1,\alpha}$``$On \cap \wq_{1}$ is cofinal in $On \cap \wq_{\alpha}$. \qed\\

%\begin{lemma} Suppose that $H(\kappa)$ is closed under sharps. Then 
%$\forall X\sset \re (X\in \Pi_{1}^{\pa {\widetilde Q,F_{\kappa}}}\equi X\in \Pi^{1}_{3})$.\end{lemma}

%\pf   \qed\\

The point of the next lemma is that although $\widetilde M$ is ostensibly about the collection of
$\Sigma_{1}$-singleton subsets of $\kappa$, with the assumption of stable
measurability, considerations about it reduce to the $\Sigma_{1}$-stable parts
of {\tmem{bounded}} subsets of $\kappa$.

\begin{lemma}
  \label{key}Suppose $\kappa$ is stably measurable. Then $\widetilde M=\widetilde Q$.
\end{lemma}

{\pf} We first remark that $\kappa$ being stably measurable implies all bounded subsets of $\kappa$ have sharps. $(\, \supseteq \,)$ is straightforward. $\left( \sset \right)$: 
$\widetilde M$ is clearly transitive.

Let  $x \in \widetilde M$ and by Lemma \ref{allaresingletons} let it be
coded by some $X \in M_{0}$.
%\todoi{Justify X} 
Let $\pa{M, \in ,F}$ witness stable
measurability. Then for some $p \in H ( \kappa )$, $\{ X \} \in
\Sigma_{1}^{M} ( \kappa ,p )$. 
Then find some $\pa{N, \in ,F_{0}} \prec \pa{M, \in ,F}$ with
$| N | < \kappa \nobracket \nobracket$, $\pa{N, \in ,F_{0}} \models$``$F_{0}$
is a normal measure on $\bar{\kappa}$'', and $X \cap \bar{\kappa} \in N_{0}$,
$p \in H ( \bar{\kappa} )^{N}$. By elementarity $\{ X \cap \bar{\kappa} \}$ is
a $\Sigma^{\pa{N, \in}}_{1} \{ \bar{\kappa} ,p \}$ singleton by the same
definition as $\{ X \}$ was. As $\pa{M, \in ,F}$ is iterable, so is $\pa{N,
\in ,F_{0}}$ and if $j_{\alpha , \beta}$ $( 0 \leq \alpha \leq \beta \in
\tmop{On} )$ are its ($\Sigma_{1}$-preserving) iteration maps, we shall have
that $\{ j_{0, \kappa} ( X \cap \bar{\kappa} ) \}$ satisfies the same
definition as that of $\{ X \}$ in $N'$ where $j_{0, \kappa} :N \imp N'$. That is: $j_{0,
\kappa} ( X \cap \bar{\kappa} ) =X$. \ Note also that $N' \in  \widetilde Q=\wmm$, as $N' \in
L_{\sigma_{N}} [ N ]$. Thus $X$ and so $x$ are in $L_{\sigma_{N}} [ N ]$ and
we are done. \ \ \ {\qed}
%\todoi{Here we don't need the good $\Sigma_{1}$-WO}

\begin{lemma}
  \label{Cor7}If $\kappa$ is stably measurable, then it is witnessed to be so
  by $(\widetilde M, \in ,F )$ where $( \widetilde M, \in )$ is as above; in particular $\la \widetilde M, \in \ra\prec_{\Sigma_{1}}\la H(\kappa^{+}),\in\ra$ itself and $F=F_{\kappa} \cap M$
  where $F_{\kappa}$ is the c.u.b. filter on $\power ( \kappa \nobracket$).
  Thus $(\widetilde M, \in ,F_{\kappa} ) \models$``$F_{\kappa}$ is the c.u.b. filter and
  is a normal measure on $\kappa$''.
\end{lemma}

{\pf} We first show that $\la \widetilde M, \in \ra\prec_{\Sigma_{1}}\la H(\kappa^{+}),\in\ra$: by assumption there is  some $\la M, \in \ra\prec_{\Sigma_{1}}\la H(\kappa^{+}),\in\ra$, some $\kappa$-model $N\supseteq M$, and some $U$ with $\la N, \in, U\ra$ witnessing stable measurability. Then $\widetilde M \sset M$ (because $M_{0}\sset M$), 
 so suppose that $\la \widetilde M, \in \ra$ is not a ${\Sigma_{1}}$ substructure of $\la M, \in \ra.$  Let $\varphi (A,\kappa, a)^{M}$ but, for a contradiction, $\neg\varphi (A,\kappa, a)^{\widetilde M}$, for some $A\sset \kappa$, $A\in \widetilde M$ where $A\in M_{0}$, and parameter $a\in H_{\kappa}$.  There is some $\psi \in \Sigma_{1}$ so that $\psi(A',\kappa,b)$ defines uniquely $A' =A$ as a $\Sigma_{1}(\kappa,b)$ singleton. By $\Sigma_1 $-elementarity, $\psi(A',\kappa,b)$ holds in $M$ and by upwards persistence both it and  $\varphi (A,\kappa, a)$ hold in $N$ too. By the same argument find $\la N',\in, U\cap N' \ra \prec  \la N, \in, U\ra$ with $\tmop{TC}(\{a\}\cup \{b\}), A \in N'\cap \kappa = \kappa_{0}\in \kappa$.  Let $\la N_{0},\in, V_{0}\ra $ be its transitive collapse with $V_{0}$ now an $N_{0}$-normal measure on $\kappa_{0}$. Then iterate  $\la N_{0},\in, V_{0}\ra $ to $\la N_{\kappa},\in, V_\kappa\ra $ with some map $j_{0,\kappa}$ now satisfying $\varphi( j_{0, \kappa} ( A \cap {\kappa}_{0} ),\kappa, a)^{N_{\kappa}}$. But $N_{\kappa}\in \widetilde M$,  and also $\psi(j_{0, \kappa} ( A \cap {\kappa}_{0} ),\kappa,b)^{N_{\kappa}}$. By uniqueness of $A$'s definition \via $\psi$ and upwards absoluteness of $\Sigma_{1}$ formulae, $j_{0, \kappa} ( A \cap {\kappa}_{0} ) = A$. But then $\varphi (A,\kappa,a)^{\widetilde M}$ - a contradiction.

We just saw that any $X \in\widetilde M \cap \power ( \kappa )$ is of the form
$j_{0, \kappa} ( X \cap {\kappa}_{0} )$ for some iteration map $j_{0, \kappa}
: ( N,F_{0} ) \imp ( N' ,F' )$ by repeating ultrapowers by an $N$-normal measure.
Thus $X=j_{0, \kappa} ( X \cap{\kappa}_{0} )$ either contains, or is disjoint
from a tail of the critical points of the embeddings $j_{\alpha , \alpha +1}$
for  $\alpha < \kappa$. As these critical points form a c.u.b
subset of $\kappa$, definable from $N \in H ( \kappa )$, and which is thus in
$\widetilde M$, $F_\kappa$ is thus a measure on $\widetilde M$. For amenability, let $\la X_{\nu}\ra _{\nu<\kappa} \in \widetilde M $ be a sequence of subsets of $\kappa$. Let it be coded by some $X\sset \kappa$, $X\in \widetilde M$, and as above have $X$ (and thus $\la X_{\nu}\ra _{\nu<\kappa}$) in some $N'$, $X=j_{0, \kappa} ( X \cap {\kappa}_{0} )$ \etc as above. $(  N' ,F' )$ is amenable and $F'$ is generated by the tail filter on the cub in $\kappa$ set of the critical points. But then $\{\nu\mid X_{\nu}\in F' \} = \{\nu\mid X_{\nu}\in F_{\kappa} \}\in N'\in \widetilde M$, and amenability is proven. The proof of $\widetilde M$-normality is similar.

Finally note that $\mb^{<\kappa}\widetilde M\sset \widetilde M$: suppose $f:\alpha \imp \widetilde M$ for some $\alpha <\kappa$. As $\widetilde M = \widetilde Q$, each $f(\xi)$ is in $L_{\sigma_{a(\xi)}}$ for some ${a(\xi)}$ a bounded subset of $\kappa$. However now code $\la a(\xi)\ra_{\xi<\alpha}$ by some $a$ still a bounded subset of $\kappa$. Then $\ran(f)\in L_{\sigma(a)}\sset \widetilde M$.
 {\qed}\\
%\todoi{Here we do need the good WO to get the $\Sigma_{1}$ elementarity}

We thus can, and do, assume that $\la \widetilde M, \inn , F_{\kappa} \cap \widetilde M\ra$  witnesses stable measurability, if it occurs.

\begin{corollary}
  $\kappa$ stably measurable implies $\la \widetilde M, \inn\ra$ is the minimal $\Sigma_{1}$-substructure of $\la H(\kappa^{+}),\in \ra$ containing  $\{\kappa\}\cup H(\kappa)$,  and $\sigma(\kappa)=On\cap \widetilde M$.
\end{corollary}

\pf Any such $\Sigma_{1}$-substructure of $\la H(\kappa^{+}),\in \ra$ must contain $\bigcup_{a \in H
  ( \kappa )} L_{\sigma_{a}} [ a ]$, which we have just seen equals $\widetilde M$. \qed

\begin{corollary}\label{sharps}
  $\kappa$ stably measurable implies that for every $A \sset \kappa $, with
$ A \in\widetilde  M $, $A^{\#}$ exists, and is in $\widetilde M$.
\end{corollary}

{\pf} Again let $A= j_{0, \kappa} ( A \cap \bar{\kappa} )$ for some iteration
$j_{0, \kappa} : ( N,F_{0} ) \imp ( N' ,F' )$. As $( N' ,F' ) \in
L_{\sigma_{N}} [ N ]$, so are the next $\omega$-many iterates $j_{\kappa ,
\kappa + \omega} : ( N' ,F' ) \imp ( \tilde{N} ,G )$ (because $( N' ,F' )\in  L_{\sigma_{N}} [ N ]$ and the latter is an admissible set); but these critical
points above $\kappa$, $\pa{\kappa_{\kappa +i} \mid 0<i< \omega}$ are Silver
indiscernibles for $L [ A ]$ and are below $\sigma_{N}$. \ Thus $A^{\#}$,
either thought of as an $A$-mouse or coded as a subset of $\kappa$, can be constructed in
$L_{\sigma_{N}} [ N ]$ and is thus in $\widetilde M$. \ \ \ \ \ \ \ \ \ {\qed}

%\todoi{Here we don't need the good $\Sigma_{1}$-WO}

\begin{lemma}\label{1.11}
%{\cor $(\forall a\sset \gamma<\kappa(a^{\sharp}\mbox{ exists }) \wedge $}
$\mbox{If there is a good } \Sigma_{1}(\kappa,p)$ wellorder of $\power(\kappa)$ for some $p\in H(\kappa)$ then:
$$ \kappa \mbox{ is stably measurable } \Equi \widetilde M = \wmm.$$
\end{lemma}
\pf The direction $(\Rightarrow)$ is Lemmata \ref{Lemma1.10} and  \ref{key} and does not require the additional assumption. For $(\Leftarrow)$:  firstly suppose that 
$\widetilde M = \wmm$; then notice trivially for 
every $a\sset \gamma<\kappa$ there is the least  $a$-mouse, $N^{a}$, witnessing that $a^{\sharp}$ exists. And its $\kappa$'th iterate $N_{\kappa}^{a}\in \widetilde Q$ ( $= \wmm$) and $\sigma_{a} < \kappa^{+ L[a]}= On \cap N^{a}_{\kappa}$. In particular $L_{\sigma_{a}}[a]\in \pa {N_{\kappa}^{a},F^{a}}$ where $F^{a}=F_{\kappa}\cap N_{\kappa}^{a}$. Consequently $\pa{\widetilde Q, F_{\kappa}}\models $``$F_{\kappa}$ is a normal ultrafilter on $\kappa$''. By the existence of the good $\Sigma_{1}$-wellorder, Lemma \ref{M} states that we have that $\widetilde M \prec_{\Sigma_{1}}H(\kappa^{+})$ and $\pa{\widetilde M, F_{\kappa}}$ winesses stable measurability.
% \\ \mb \
%Secondly if we assume that some bounded subset $a$ of $\kappa$ has no sharp, then $a\notin \wq$ but $a$ is in $\widetilde M$. So the right hand side fails and the implication is trivially true.
%If we assume that $\kappa$ is not stably measurable, then we must have that $\widetilde M \neq \widetilde Q.$ Secondly, suppose some $a\sset \gamma <\kappa$ has no $\sharp$. But then $a\in L_{\sigma_{a}}[a]\sset \widetilde M\backslash \widetilde Q$.
\qed\\

%\todoi{Check the boldface version of good wellorders in the last lemma}

In fact there is more to be said on the sharps in $\widetilde M$. 
%By the last corollary if $\kappa$ is stably measurable the supposition of the next definition is fulfilled. 
\begin{lemma}
  \label{u=sigma}Let $\kappa$ be stably measurable. Then $u_{2} ( \kappa ) =
  \sigma ( \kappa )$.
\end{lemma}

\nod{\bf Proof:}
  %Let $(\widetilde M, \in ,F_{\kappa} )$ witness stable measurability. 
  $( \leq )$ Let $a \in H (
  \kappa )$ be a set of ordinals. Then  $a^{\#}$ (which exists by Cor. \ref{sharps}), considered as the least
  $a$-mouse $( \bar{N}^{a} , \in , \bar{U} )$ is in $H ( \kappa \nobracket$)
  and can be iterated $\kappa +1$ \ many times, inside
  $L_{\sigma_{\alpha^{\#}}} [ a^{\#} ] \sset \widetilde M$. If these iterations points
  are $\{ \lambda_{\alpha} \}_{\alpha \leq \kappa +1}$ then as above these are
  Silver indiscernibles for $L [ a ]$ and thus $\lambda_{\kappa +1} =
  \iota_{\kappa +1}^{a} < \sigma_{\alpha^{\#}} < \tmop{On} \cap \widetilde M= \sigma .$
  
  $( \geq )$ Just note that for any $\gamma < \sigma = \tmop{On} \cap\, \widetilde M$ there is, by Lemma
  \ref{key}, some $a \in H ( \kappa )$ with $\gamma < \sigma_{a} \leq \sigma .$
  But $a^{\#}$ exists and then $\gamma < \sigma_{a} <
  \kappa^{+L [ a ]} <  u_{2}$.%}\sigma_{a^{\#}} < \sigma$. 
\qed\\

However the converse of the last lemma may fail: suppose ($\kappa = \omega_1$) that $u_2(\omega_{1}=\omega_2$ (which it may, by a result of Woodin, if there is a measurable cardinal and $NS_{\omega_1}$ is saturated); but then also $\sigma(\omega_1) = \omega_2$. It is easy to see that $\kappa$ stably measurable implies that $\kappa$ is Mahlo.
Hence in general $u_2(\kappa) = \sigma(\kappa) \not \imp \kappa$ is stably measurable.\\

%\todoi{Q. Can this fail for inaccessibles rather than just for $\kappa = \omega_1$?}
%\begin{lemma}
%\end{lemma}

%\todoi{Here we use the good $\Sigma_{1}$-WO. to get $\widetilde M$ to be stable, and so $\sigma = On \cap \widetilde M $ for the $(\geq)$ direction only.?  Mmmm Think not:  $\sigma = On \cap \widetilde M $ anyway. }
%\rem One might hope for a converse to the last lemma assuming all bounded subsets of $\kappa$ have sharps,  and there is a good $\Sigma_1$-wellorder of $\power(\kappa)$, but although we have $\sigma = On \cap \widetilde M$ and $u_{2}=On \cap \widetilde Q$ from this we do not quite seem to  get in general $u_{2} < \sigma$ from the failure of stable measurability. 

The following is similar to L\"ucke 7.1(ii) showing weakly compact cardinals with the {\boldmath{$\Sigma_{1}$}}-club property (to be defined below) reflect on a stationary set.
 
\begin{lemma} If $\kappa$ is weakly compact and stably measurable, then the set of cardinals $\alpha$ below $\kappa$ which are stably measurable is stationary.
\end{lemma}
\pf Let $\la \widebar  M^{\kappa},F_{\kappa}\ra$ witness the stable measurability of $\kappa$. Thus $\widebar  M^{\kappa} \prec _{\Sigma_{1}} H(\kappa^{+})$. Let $C\sset \kappa$ be cub.  Choose $M \prec  H(\kappa^{+})$ with $|M|=\kappa$ and $ \widebar  M^{\kappa}\cup \{ \widebar  M^{\kappa}, C\}\sset M$ and $ \mb ^{<\kappa}M \sset M$ with some elementary map $j:M\imp N$, with critical point $\kappa$ as given by weak compactness. Note that $ \widebar  M^{\kappa} \prec _{\Sigma_{1}} M$. In general $H(\kappa^{+}))^{M}\subsetneq ( H(\kappa^{+}))^{N}$, but $\power(\kappa)^{M}\sset N$ (and $(F_\kappa)^M\sset (F_\kappa)^N$). As $\widebar  M^{\kappa}$ is an element of $H(\kappa^{+}))^{M} $
%of cardinality $\kappa$ in $M$, it can be coded by a subset $E\sset \kappa$; $E$ and 
it is in $N$. We claim:\\

\nod {\em Claim:} {\em  $ \widebar  M^{\kappa} \prec _{\Sigma_{1}} ( H(\kappa^{+}))^{N}$ and thus $ \la \widebar  M^{\kappa},F_{\kappa}\ra$ witnesses stable measurability of $\kappa$ in $N$.}

If the claim holds: $$N\models j(C)\cap \{\alpha <j(\kappa)\mid \exists 
\widebar  M^{\alpha} \prec _{\Sigma_{1}} ( H(\alpha^{+})), 
\la \widebar  M^{\alpha},F_{\alpha}\ra \mbox{ witnesses stable measurability } \} \neq \emp.$$
But then there is some $\alpha \in C$ with $ \la \widebar  M^{\alpha},F_{\alpha}\ra$ witnessing stable measurability, and we are done. \\
Proof of Claim: Let $\vec A\in \widebar  M^{\kappa}, \,\varphi \in \Sigma_{1}$ with $\varphi(\vec A)^{N}$. By upwards absoluteness: $\varphi(\vec A)^{H(\kappa^{+})}$ and then by downwards $\Sigma_{1}$-elementarity: $\varphi(\vec A)^{ \widebar  M^{\kappa}}$.  \qed (Claim \& Lemma)\\

The next result says that stable measurability is easily propagated upwards; but is perhaps less surprising when one realises that stable measurability at $\kappa$ is more about the bounded subsets of $\kappa$. \cite{Luecke18} Thm. 7.4 has that  a stationary limit of iterable cardinals has the 
 {\boldmath${\Sigma}  _{1}$}-$\tmop{club}$ property (to be defined below). We have a weaker hypothesis and a stronger conclusion.

\begin{theorem}If $\kappa$ is the stationary limit of stably measurable cardinals, then 
 $\kappa$ is stably measurable.
\end{theorem}

\pf 
% We first suppose there is  wellorder of $\power(\kappa)$ definable over
% $H(\kappa^{+})$ in order to define Skolem functions, or just directly,
%a $\Sigma_{1}$-satisfaction predicate $S$ for $H(\kappa^{+})$. We may assume this is $\Sigma_{n}^{{H(\kappa^{+}})}(\vec z)$ for some 
Using $AC$, choose $S$ a  $\Sigma_{1}$-satisfaction predicate  for $\la H(\kappa^{+}),\in \ra$.
%For some $N\in \omega$ then we have we may choose 
Choose $\pa {X,\in, S\cap X}\prec \pa {H(\kappa^{+}),\in, S}$ with {\cor  $\vec z $}, $X\cap \kappa \in \kappa$, and $H(X\cap \kappa)\sset X$ (note $\kappa$ is a strong limit), and letting $\pi: \pa {X, X\cap S} \imp \pa {\bar H, \bar S}$ be the transitive collapse, let $\pi(\kappa) = \bar \kappa$. By assumption we may additionally assume that $\bar \kappa $ is stably measurable. Then, let $\bar M =  \bigcup_{a \in H
  (\bar \kappa )} L_{\sigma_{a}} [ a ]= \wmm(\bar \kappa) =\wq(\bar \kappa)$ (the latter since by assumption all bounded subsets of $\kappa$ will have sharps); the sets of the right hand side here are all contained in $\bar H$. Then $\pa {\bar M, F_{\bar \kappa}\cap \bar M}\in \bar H$ and is definable there. By the stable measurability of $\bar \kappa$, \ie using that $\bar M \prec_{\Sigma_{1}}H(\bar \kappa^{+})$, and the inclusion $\bar M \sset \bar H \sset H(\bar \kappa ^{+})$, and noting that $\bar S$ codes $\Sigma_{1}$-satisfaction over $\pa{\bar H,\in}$, we have that 
  $$\pa{ \bar H,\bar S} \models \, \bar M \prec_{\Sigma_{1}} V\, \wedge \pa{\bar M, F_{\bar\kappa}}\models\mbox{``} F_{\bar\kappa}\mbox{ is a normal measure on }\bar\kappa\mbox{ ''.}$$ 
 % But this is a $\Sigma_{N}^{\bar H}(\vec z,\bar \kappa)$-property of $\pa {\bar M, F_{\bar \kappa}}$. 
 Applying $\pi^{-1}$ we have  $\pi^{-1}(\pa {\bar M, F_{\bar \kappa}}) = \pa {\wq(\kappa), F_{ \kappa}}$.
 % where as above $\widetilde M = \bigcup_{a \in H
 % ( \kappa )} L_{\sigma_{a}} [ a ]$. 
 We then have:
   $$  \pa {H(\kappa^{+}),S}\models \, \wq(\kappa) \prec_{\Sigma_{1}} V\, \wedge \pa{\wq(\kappa), F_{\kappa}}\models\mbox{``} F_{\kappa}\mbox{ is a normal measure on }\kappa\mbox{ ''.}$$ 
In other words, $\pa{\wq(\kappa), F_{\kappa}}$ witnesses that $\kappa$ is stably measurable.\qed\\

We now relate stable measurability and its analysis above to L\"ucke's notion of  the {\boldmath${\Sigma}  _{1}$}-$\tmop{club}$ property.

\begin{definition}
  (L\"ucke {\cite{Luecke18}} Lemma 4.1)\label{sigmaclub} $\kappa$ has the
  {\boldmath${\Sigma}  _{1}$}-$\tmop{club}$ property if, for any $A \sset \kappa$
  so that $\{ A \} \in \Sigma^{}_{1} ( \kappa ,z )$ where $z \in H ( \kappa
  )$, then $A$ contains or is disjoint from a club subset of $\kappa$. 
\end{definition}

(Actually this is not L\"ucke's basic definition, but he shows this is
equivalent to it.) \ Note that by `$\Sigma_{1} ( \kappa ,z )$ definable', we
can take this to mean $\Sigma_{1}^{H ( \kappa^{_{+}} )} ( \kappa ,z
)$-definable, by L\"owenheim-Skolem and upwards absoluteness arguments.

We introduced in {\cite{Sharpe-Welch2012}} the following notion when
discussing variants of Ramseyness.

\begin{definition}\label{Sigmaclub}
  $\kappa$ is called {\tmem{($\omega_{1}$-)iterable}} if for any $A \sset
  \kappa$ there is a transitive set $M$, and filter $U$, with $A \in M$ and $(
  M, \in ,U ) \models $``$ U$ is a normal measure''; it is amenable, iterable by $U$
  and has wellfounded ultrapowers.
\end{definition}

(In {\cite{Sharpe-Welch2012}} this was rather obscurely called the $Q$ property.) It was shown
there ({\em op. cit.} Lemma 5.2) to be strictly weaker than Ramseyness: that would require additionally
that the filters $U$ be $\omega$-closed. One can show that an
$\omega_{1}$-Erdos cardinal is a stationary limit of $\omega_{1}$-iterable
cardinals (see {\cite{Sharpe-Welch2012}} Lemma 5.2). But notice that iterability is
clearly stronger than stable measurability: {\tmem{every}} subset of $\kappa$
must be in some iterable structure, not just the $\Sigma_{1} ( \kappa
)$-singletons.

L\"ucke shows the following:

\begin{theorem}
  (L\"ucke {\cite{Luecke18}} Cors. 4.12 and 4.5)  (i) $\kappa$ iterable
  $\Imp  $ the $\tmmathbf{\Sigma}_{1}$-$\tmop{club}$ property holds at
  $\kappa$.
  
  (ii) The $\tmmathbf{\Sigma}_{1}$-$\tmop{club}$ property at $\kappa
  \hspace{1em}   \Imp   \hspace{1em} \all x \in \mathbbm{R} ( x^{\#}
  \nobracket$ exists).
\end{theorem}

We remark later that the gap above can be closed by showing that the
$\tmmathbf{\Sigma}_{1}$-$\tmop{club}$ property is equiconsistent with stable
measurability. However first we may show outright:

\begin{theorem}
  \label{main}$\kappa$ has the $\tmmathbf{\Sigma}_{1}$-$\tmop{club}$ property,
   if $\kappa$ is stably measurable.
\end{theorem}

{\pf}%\todoi{This now looks wrong. Only R-to-Left looks OK. For L-to-R looks OK?? if we add there is a $\Sigma_{1}^{H(\kappa^{+})}(\kappa)$ good wellorder of $\power(\kappa)$} 
%$\left( \Imp \right)\widetilde  M \cap \power ( \kappa )$ consists only of
%$\Sigma_{1}(\kappa, p)$-singletons $\{ A \}$.
% of the form of the last definition. 
%Hence $( \widetilde M, \in ,F_{\kappa} \cap \widetilde M )$ witnesses stable measurability.
%$( \Leftarrow )$ 
We've seen above at Corollary \ref{Cor7} that if $\kappa$ is
stably measurable, then it is witnessed to be so by $(\widetilde M, \in ,F_{\kappa} \cap
\widetilde M )$; but the latter contains $ M_{0}$ so this suffices. \ {\qed}\\

The converse can be false:

\begin{lemma}\label{L1.27} $ZFC\not \vdash \kappa $ has the \siclub \, property $\imp \kappa $ is stably measurable. 
\end{lemma}
\pf L\"ucke points out in \cite{Luecke18} Cor. 7.3, that if $\kappa$ is a regular limit of measurables, then the $\tmmathbf{\Sigma}_{1}$-$\tmop{club}$ property holds. But such a $\kappa$ need not be Mahlo, and so not stably measurable. \qed\\

Conversely we now have (and by the above  the assumption in the lemma is necessary):
\begin{lemma}\label{1.21}
Assume there is a good {\boldmath$\Sigma_{1}$}$(\kappa)$-WO of $\power(\kappa)$. 
 Then $\kappa$ has the $\tmmathbf{\Sigma}_{1}$-$\tmop{club}$ property
   implies $\kappa$ is stably measurable.%, as witnessed by $\pa{\widetilde M,\in, F_{\kappa}}$.
\end{lemma}

\pf  That $\kappa$ has the $\tmmathbf{\Sigma}_{1}$-$\tmop{club}$ property ensures, by an application of Lemma \ref{allaresingletons} that $F_{\kappa}$ measures $\power(\kappa)\cap \widetilde M$. That there is a good {\boldmath$\Sigma_{1}$}$(\kappa)$-WO of $\power(\kappa)$ will ensure that $\widetilde M \prec_{\Sigma_{1}}H(\kappa^{+})$.  \qed\\

Putting the argument of the last lemma together wth the fact that stably measurable cardinals are Mahlo, one can conclude:

\begin{corollary} If $\kappa$ is a regular cardinal which is not Mahlo, but is limit of measurable cardinals, then there fails to be a good {\boldmath$\Sigma_1$}$(\kappa)$-wellorder of $\power(\kappa)$.
\end{corollary}

In fact \cite{LS18} Cor. 1.4 show this directly for lightface $\Sigma_1(\kappa)$ good wellorders, but for all regular limits of measurables.

\section{Stable measurability in $L[E]$-models}

We consider what happens when stable measurability is instantiated in models with fine structure. The outcome is an equivalence between the notions considered.

\subsection{When $K=K^{DJ}$}

We let in this subsection $K=K^{\tmop{DJ}}$. 
%L\"ucke showed
%\begin{theorem}
 % ({\tmem{{\cite{Luecke18}}}} 4.10) that the $\tmmathbf{\Sigma}_{1}$-$\tmop{club}$ property is downward absolute to $K.$
%\end{theorem}
We shall show that the stable measurability is  downward absolute to $K$. 

We note first:
\begin{lemma}\label{1.28}$(V=K^{DJ})$
$\widetilde M \prec_{\Sigma_{1}}H(\kappa^{+})$. 
\end{lemma}
\pf By Lemma \ref{M}, as in $K^{DJ}$ we have a good $\Sigma_{1}^{H(\kappa^{+})}(\kappa)$ wellorder $<$ of $\power(\kappa)$.
\qed\\

We then relate $\wq (\kappa)$ to an older notion.
% indeed the forbear of $\wb(\kappa)$.% of `$Q$-structure at $\kappa$'.

%Hence (or it can be shown directly - see the proof of the next Theorem):

%\begin{corollary} Stable measurability is downward absolute to $K$.
%\end{corollary}

\begin{definition}
  (The $Q$-structure at $\kappa$)(\cite{D}) In K, let $Q ( \kappa ) \dfs
  \pa{J^{F_{\kappa}}_{\theta ( \kappa )} , \in ,F_{\kappa}}$ be the union of
  the $\kappa$'th iterates of all DJ-mice $M \in H ( \kappa )$.
\end{definition}

As the measure of each such $\kappa$-iterate $M_{\kappa}$ of such a DJ-mouse $M \in H
( \kappa )$, is generated by the tail sequence filter of its closed and
unbounded in $\kappa$ sequence of critical points, the measure on $ M_{\kappa}$
is just $F_{\kappa} \cap M_{\kappa}$, and thus $M_{\kappa}$ is of the form
$\pa{J^{F_{\kappa}}_{\alpha} , \in ,F_{\kappa}}$. $Q ( \kappa )$ is the union
of all such, and is itself a DJ-mouse. (The reader should be reminded that DJ-mice, whilst amenable, are not acceptable in the modern meaning of the word. Indeed for a DJ-mouse $M$ with critical point $\kappa$ it need not be the case that $(H_\kappa)^M\in M$. Such is the case for example with $Q(\kappa)$.) The height of $Q ( \kappa )$ is thus
proportional to the length of the critical mouse order of $H ( \kappa )$. \ \
(It can be shown (i) that  $\tmop{if}   \eta$ is this order type then $\theta (
\kappa ) = \kappa \cdot \eta$, and thus (ii)  $H(\kappa)$ is closed under sharps iff $\eta$ is a multiple of $\kappa^{2}$.) 

\begin{lemma}\label{utheta}
In $K^{DJ}$: for any cardinal $\kappa$, $\theta(\kappa)\leq u_{2}(\kappa)$.
\end{lemma}

Still in $K^{DJ}$, \cite{W3} Lemma 3(i) shows that the uniform indiscernibles for bounded subsets of $\kappa$ (of which thus $u_{2}(\kappa)$ is the second) are precisely the critical points of the successive ultrapowers of $Q(\kappa)$.  $Q(\kappa)$ need not have the all the sets of $\wq(\kappa)$ (it may be too short, indeed in this case even if all bounded subsets of $\kappa$ have sharps, we may have $\wq(\kappa)\neq Q(\kappa)$) but if $Q(\kappa)$ is admissible then we shall have $Q(\kappa) = \wq (\kappa)$.  Still assuming $Q(\kappa)$ is admissible the discussion in \cite{W12} showed that $u_{2}(\kappa)=\theta(\kappa)$.
What we shall see is that if in $K$, $\sigma (
\kappa ) =u_{2} ( \kappa )$, then we shall have also that $\theta ( \kappa ) =
\sigma ( \kappa )$ and moreover that $\wq(\kappa)=Q ( \kappa ) =
\pa{J^{F_{\kappa}}_{\theta ( \kappa )} , \in ,F_{\kappa}}$ itself witnesses
stable measurability in $K$. 

\begin{lemma}
Suppose $V= K^{DJ}$ and that $Q(\kappa)$ is admissible. Then $\wq (\kappa) = Q(\kappa)$.
\end{lemma}
\pf It is easy to see that $(\,\supseteq\,)$ holds, by the previous style of arguments. For $(\sset)$: let $a \in H(\kappa)\cap \power(\kappa)$. The $a$ is simply an element of a $DJ$-mouse $N\in H(\kappa)$ (as $K^{DJ}$ is the union of such). However then $a\in N_{\kappa}$ which is an initial segment of $Q(\kappa)$. Now suppose $x\in \wq $; then $x \in L_{\sigma_{a}}[a]$ for such an $a$. (We are using here, that as $Q(\kappa)$ is admissible, $On\cap Q(\kappa)$ is a multiple of $\kappa^{2}$ and thus $H(\kappa)$ is certainly closed under $\sharp$'s, and thus $\wmm = \wq$.) Then $a^{\sharp}$ is in some $DJ$-mouse $M\in H(\kappa)$. But $Q(\kappa)\supseteq H(\kappa)$. Hence $M,\kappa \in Q(\kappa)$. By KP then the $\kappa$'th iterate of $M$, $M_{\kappa}$ is  in $Q$. But $\power(\kappa)\cap L[a]\sset M_{\kappa}$. Thus there is a subset of $\kappa$ that codes the ordinal $\sigma_{a}$, and so also a code for the structure $L_{\sigma_{a}}[a]$, in $M_{\kappa}$, and so, by KP again, these sets themselves are in $Q$.  This puts $x\in Q$.
\\ \mb \hfill\qed \\
%We also have:

\begin{theorem}\label{ordinals}
  (i) $\sigma ( \kappa ) =u_{2} ( \kappa )   \Imp$ $\sigma ( \kappa )^{K}
  =u_{2} ( \kappa )^{K} $. If additionally $\neg 0^{\dagger}$ \ then
  $\sigma ( \kappa ) = \sigma ( \kappa )^{K}$.\\
  (ii)  $\kappa \nobracket$ is stably
  measurable $ \Imp   ( \kappa \nobracket$ is stably
  measurable$\nobracket )^{K}$  as witnessed by $Q ( \kappa ) =
\pa{J^{F_{\kappa}}_{\theta ( \kappa )} , \in ,F_{\kappa}}$.
\end{theorem}

{\nod}{\pf}For (i):  assume $\sigma ( \kappa ) =u_{2} ( \kappa )$. Firstly note that if
$0^{\dagger}$ exists, then every uncountable cardinal $\kappa$ is Ramsey in
$K$, and hence is iterable, hence stably measurable in $K.$ Then the
conclusion follows by Lemma \ref{u=sigma}. \ So assume $\neg 0^{\dagger}$.

(1) $\sigma ( \kappa )^{K} = \sigma ( \kappa ) =u_{2} ( \kappa ) =u_{2} (
\kappa )^{K} .$

{\pf}of (1). By $\Sigma^{1}_{3}$-absoluteness arguments going back to Jensen (see, \eg, \cite{D} or \cite{DJK})
for reals, but applying them for bounded subsets of $\kappa$, $u_{2} ( \kappa ) =u_{2} ( \kappa )^{K}$. So we are left with
showing the following Claim:

{\tmem{Claim }}$\sigma ( \kappa )^{K} = \sigma ( \kappa )$

{\pf} $\sigma ( \kappa )^{K} \leq \sigma ( \kappa )$ follows from the
wellorder of $\power ( \kappa ) \cap K $ being a good $\Sigma^{H ( \kappa^{+}
)^{K}}_{1} ( \kappa )$-definable wellorder which at the same time is a good
$\Sigma^{H ( \kappa^{+} )}_{1} ( \kappa )$ wellorder in $V$; thus if $\{ A \}$ is a
$\Sigma_{1} ( \kappa ,p )^{K}$ singleton subset of $\kappa$, it is also a
$\Sigma_{1} ( \kappa ,p )^{}$ singleton in $V$. Hence any such $A \in  M_{0}^{K}$
coding a wellorder $\tau < \sigma ( \kappa )^{{K}}$ is also in $ M_{0}$. Clearly then
$\tau$ and so $\sigma ( \kappa )^{K} \leq \sigma ( \kappa )$.

But $\sigma ( \kappa )^{K} \geq u_{2} ( \kappa )^{K}$, since the latter is also
$\sup \{ \tmop{cp} ( N_{\kappa +1} ) \mid  N_{\kappa +1} \nobracket$ is the
$\kappa +1$'st iterate of a mouse $N$ in $ H ( \kappa ) \}$ and moreover
$On\cap N_{{\kappa +1}}< \sigma_{N}$.
All
such $N_{\kappa +1}$ are in $\bar M$ if the latter is any
$\Sigma_{1}$-substructure of $H ( \kappa^{+} )^{K}$ containing $H ( \kappa )
\cup \{ \kappa \}$. Hence $\sigma ( \kappa )^{K} \geq  u_{2} ( \kappa )^{K} = u_{2}(\kappa)=\sigma ( \kappa )$. 
\\\mb \hfill {\qedtwo{Claim \& (i)}}\\

\nod For (ii) assume that $\kappa \nobracket$ is stably
  measurable.
  
{\tmem{Claim }} $Q ( \kappa \nobracket$) witnesses that $\kappa$ is stably
measurable in $K$.

{\pf} Work in $K$. Let $\widetilde M=\widetilde M^{K}$. \ $Q ( \kappa ) \sset \widetilde M$ \ since $Q ( \kappa
)$ is the union of the $\kappa$'th iterate of DJ-mice $N \in H ( \kappa )$ and
all such iterates are in $\widetilde  M$.

$Q ( \kappa ) \supseteq \widetilde  M$: By Lemma \ref{allaresingletons} it suffices to show $M_{0}^{K}\sset Q(\kappa)$. Let $A \in M_{0}^{K}$. 
By the argument for (1), $A \in M_{0}$, and by Corollary \ref{sharps} , using stable measurability in $V$,  $A^{\#}$ exists, and by absoluteness it exists in $K$.

Hence $A \in M_{0}^{K} \cap \power ( \kappa )\, \Imp \,A^{\#} \in M_{0}^{K}$. However then there is
some DJ-mouse $N_{A}$ with $A \in N_{A}$. Note now the $<^{\ast}$-least such mouse $N_{A}$ projects to $\kappa$ and so has a code $B$ a subset of $\kappa$. But $\{A\}$ is a $\Sigma_{1}(\kappa,p)$ singleton set (some $p\in H(\kappa)$), and thus such a code set $\{B\}$ is also a $\Sigma_{1}(\kappa,p)$ singleton set and so it, and thence $N_{A}$, is in $\widetilde M$.

%(by $\widetilde M\prec_{\Sigma_{1}} H ( \kappa^{+} )$).
%\todoi{Warning: cannot appeal to this.} 
Moreover if $\pa{\lambda_{i} \mid i \in
\omega}$ are the first $\omega$ iteration points of $N_{A}$ which are Silver
indiscernibles for $L [ A ]$, then $\tilde{\lambda} = \sup \{ \lambda_{i}
\}_{i< \omega} < \sigma =u_{2} ( \kappa )$ (the latter equality by part (i)). So there is some $\bar{N} \in H (
\kappa )$ with $\tmop{cp} ( \bar{N}_{\kappa +1} ) > \tilde{\lambda}$. As
$\bar{N}_{\kappa +1}$ is a DJ-mouse, there is some $f: \kappa \imp \tmop{On}
\cap \bar{N}_{\kappa +1}$ which collapses $\tilde \lambda$ with $f \in
\Sigma_{\omega} ( \bar{N}_{\kappa +1} )$. In particular $\tilde{\lambda}$ is
collapsed, so $\bar{N} \mbox{ } ^{\ast}\! \!\geq N_{A}$. However then $A \in \power ( \kappa
)^{N_{A}} \sset \power ( \kappa )^{\bar{N}_{\kappa +1}} \sset \power ( \kappa
)^{Q ( \kappa )}$. Thus $Q(\kappa) = \widetilde M$ and $\la Q(\kappa), F_\kappa\ra$ is iterable \etc  So $\kappa$ is stably measurable.  \mb \hfill {\qedtwo{Claim \& (ii) \& Theorem}}\\

%This completes (ii). \qed

\begin{theorem}\label{ordinalsinK}$(V=K^{DJ})$ 
$$\sigma ( \kappa ) =u_{2} ( \kappa )   \Equi \kappa\mbox{ has the
  {\boldmath{${\Sigma}_{1}$}}-}\tmop{club}\mbox{ property }\Equi \kappa \mbox{ is stably
  measurable. }$$
\end{theorem}
\pf
 Note first that $\widetilde M \prec_{\Sigma_{1}}H(\kappa^{+})$. This is by Lemma \ref{M} as in $K^{DJ}$ we have a good $\Sigma_{1}^{H(\kappa^{+})}(\kappa)$ wellorder  of $\power(\kappa)$. %This means we can take the $<$-least witness to any $\Sigma_{1}$ sentence about any $A\in \widetilde M$ as itself coded by a subset $B\sset \kappa$ with $\{B\}$ a $\Sigma_{1}$-singleton. That, is with $B\in \widetilde M$. 
 
 If $\kappa$ has the {\boldmath{${\Sigma}_{1}$}}-$\tmop{club}$ property then $\la \widetilde M, \in, F_{\kappa}\ra \models$``$F$ is a normal measure on $\kappa$'', and as usual is iterable.
  Thus  $\la \widetilde M, \in, F_{\kappa}\ra$ witnesses stable measurability. This in turn implies $\sigma =u_{2} ( \kappa )$ (by \ref{u=sigma}). We are left with showing $\sigma =u_{2} ( \kappa )$ implies  the {\boldmath{${\Sigma}_{1}$}}-$\tmop{club}$ property.  As we have  $\widetilde M \prec_{\Sigma_{1}}H(\kappa^{+})$, it suffices to show that $F_{\kappa}$ measures all $\power(\kappa)^{\widetilde M}$. 
%%%%%%%%%%

Let $A \in \power(\kappa)^{\widetilde M}$.
Then $A \in M \Imp \sigma_{A}
\leq \sigma$. As we are in $K$ if $\neg A^{\#}$, then $K=L [ A ]$. (If we
define $K^{L[A]}$ inside $L [ A ]$ and this is not all of $K$, then there is some
least mouse $P \notin L [ A ]$. But then $P$ generates $A^{\#}$.) \ But in
this case, as $H ( \kappa ) =H ( \kappa )^{L [ A ]}$ we should have that if
$\omicron(A)$ is the least ordinal so that $\mathcal{A=} L_{\omicron(A)} [ A ]
\models \tmop{KP}$, $o \geq u_{2}$, as all $\kappa+1$'st iterates of mice $N \in H
( \kappa )$ are in fact in $\mathcal{A} \nosymbol$. But $\mathcal{A}$ is
merely the first $A$-admissible $> \kappa$ containing $H ( \kappa ) \cup \{
\kappa \}$. Thus $\omicron(A) < \sigma_{A}$ (as $\sigma_A$ is a limit of $A$-admissibles) and the latter is $\leq \sigma = u_{2}$ - a
contradiction. Hence $A^{\#}$ exists. Let $N_{A}$ be the $<^{\ast}$-least mouse with $A \in N_{A}$. By the $\Sigma_{1}$ elementarity of $\widetilde M$, we have $N_{A}\in \widetilde M$. By the same argument with $N_{A}$ in place of $A$ we cannot have $H_{\kappa}\sset L_{\omicron(N_A)}[N_{A}] $ the least admissible set containing $N_{A}$. Hence there is some $<^{\ast}$-least mouse $\bar M\in H_{\kappa}\backslash L_{\omicron(N_A)}[N_{A}]$. Thus $N_{A}<^{\ast} \bar M$. As $A\in  \power(\kappa)^{N_{A}} \sset
\power(\kappa)^{\bar M_{\kappa}}$ where ${\bar M_{\kappa}}$ is the $\kappa$'th iterate of $\bar M$, either $A$ or $cA$ contains a tail of the club of critical points $C_{\bar M}\sset \kappa$.\hfill\qed \\

%%%%%%%%%%

%this follows from (ii) by Lemma \ref{u=sigma} and Theorem \ref{main}. \ 

\begin{corollary}\label{2.7}
  In $K^{DJ}$, if $\sigma ( \kappa ) =u_{2} ( \kappa )$ then these two ordinals
  both equal $\theta ( \kappa ) \nocomma$ and if $<^{\ast}$ is the
  prewellordering of mice, then $o.t. \left( <^{\ast} \rest H_{\kappa} \right)
  =  \sigma ( \kappa \nobracket$). 
\end{corollary}

{\rem}In $K^{DJ}$ it can happen that $\theta ( \kappa ) < \sigma ( \kappa
\nobracket$) (for example if $K=L [ 0^{\#} ]$) but $\theta ( \kappa )$ can
never be strictly greater than $\sigma ( \kappa \nobracket$) as we always have
$Q ( \kappa ) \sset \widetilde  M_{0}$.
%\todoi{Looks like this should be $\widetilde M$}
Now just as a corollary to the above we have immediately:

\begin{theorem} The following are equiconsistent over $ZFC$:\\
(i)\,\, \, $\ex \kappa( \kappa \mbox{ is stably measurable})$ ;\\
(ii)\,\,\,
$\ex \kappa( ${\boldmath{$\Sigma_{1}$}}-club property holds at $\kappa)$ ;\\
(iii)\, \,$\ex \kappa( \sigma(\kappa) = u_{2}(\kappa)).$
\end{theorem}

Philipp L\"ucke has also pointed out that a further equivalence  can now be obtained in $K^{DJ}$ with a hypothesis that is also used in his paper \cite{Luecke18} at Lemma 3.13 and Theorem 3.14. We derive this as follows.

\begin{lemma} In $K^{DJ}$ we have $\kappa$ is  stably measurable iff $H(\kappa)$ is  not {\boldmath{$\Sigma_{1}$}}$(\kappa)$-definable.
\end{lemma}

\pf Note that $H(\kappa)\sset Q(\kappa)$ and is a $\Sigma_{1}$-definable class over, but is never an element of, the latter. By definition of $\widetilde M$ we always have $Q(\kappa)\sset \widetilde M$.  Hence the equivalences $On\,\cap\, Q(\kappa) < \sigma(\kappa)=On\cap \widetilde M$ iff $Q(\kappa)\in \widetilde M$ iff $H(\kappa)\in \widetilde M$  iff $H(\kappa)$ is {\boldmath{$\Sigma_{1}$}}$(\kappa)$-definable
 are  all true for any $\kappa > \omega$.
 However by Theorem \ref{ordinalsinK} and Cor. \ref{2.7} we have 
 $On\,\cap\, Q(\kappa) = u_{2}(\kappa)=\sigma(\kappa)$ iff $\kappa$ is stably measurable. \qed\\

\subsection{When $K=K^{{strong}}$}

In this subsection we assume $V=K$ but $\neg 0^{pistol}$. There is thus no mouse $M$ with a measure with a critical point $\kappa$ and $\lambda < \kappa$ with $o^M(\lambda)\geq \kappa$.  (Such a mouse engenders a sharp for an inner model with a strong cardinal.) Let us call $K$ built under this hypothesis $K^{strong}$.

%\bu\, Looks like in $K$: $u_{2}=\sigma(\kappa) \Imp \kappa $ stably meas. as before.\\

%Q What about $u_{2}=\sigma(\kappa))^{V}$? Difficulty is to show $\sigma(\kappa)^{K}\leq \sigma(\kappa)$.

\begin{lemma}\label{strongwo}
Suppose the measurable cardinals in $K$ 
%satisfies the following property, with $E=E^K$ 
are bounded by some $\lambda^{+} <\kappa$.
% $$(\ast)\quad \exists \lambda^{+} < \kappa \forall \nu \in (\lambda^{+} ,\kappa)( \tmop{crit}(E_\nu) \mbox{  is not a measurable cardinal} > \lambda^{+}).$$

%neither a limit of measurable cardinals, nor is `overlapped' by an extender on the $E^K$-sequence (that is, there is no $\nu >\kappa>\lambda$ with $E^K_\nu$ indexing an extender with critical point $\lambda$). 

Then there is a good {\boldmath$\Sigma_1$}$(\kappa)$-wellorder of $\power(\kappa)$. 
\end{lemma}

Note the assumption here implies that although the measurable cardinals of $K$ below $\kappa$ are bounded by some $\lambda^{+}$,
but allows some measurable $\tau\leq\lambda$ to be strong up to $\kappa$.\\

\pf Let $e\dfs E^K\rest \lambda^{+}$ with $\lambda$ such a bound. If some $\tau \leq \lambda$ is strong up to $\kappa$ on the sequence $E^{K}$, then, by a use of $\neg 0^{pistol}$ we may take $\lambda$ as this $\tau$. Then let $\psi(E,\lambda)$ be the assertion that $\lambda$ is strong up to $\kappa$ as witnessed by the sequence $E$. Otherwise let $\psi(E,\lambda)$ be ``All measurable cardinals on the sequence $E$ have their critical points $\leq \lambda$''. Then $e$ will serve as a parameter for defining the wellorder on $\power(\kappa)$ given by: \\
%\begin{array}{cl}

 We shall set   $ x \triangleleft y \mbox{ iff } x<_M y$
%  \end{array}
where $<_M$ is the usual order of construction of the structure of $M$, for an $M$ satisfying the following: 

``$E^M\rest \lambda^{+} = e\wedge
M\models KP + 
%(\ast)+ 
\psi(E^{M},\lambda)\, \wedge $ 
$  \mbox{ \em $M$ is a sound mouse } \wedge $ \\ 
$  M\, \mbox{\em  is the least level of the $J^{E^M}$-hierarchy that contains $x$ and $y$ which is a $KP$-model }\, \wedge \, $
  %\rho^\omega_M=\kappa 

   $\wedge \, M$ {\em  is $\in$-minimal satisfying these conditions}''.\\
   
   Note that these conditions require that $\rho^\omega_M=\kappa $.
% Suppose first that the measurables of $K$ are bounded by the cardinal $\lambda <\kappa$. Set $\widebar E= E^K\rest\lambda$. Then in $K$ we may build a $K$-sequence over the predicate $\widebar E$ of $\widebar E$-mice. Then we shall have that $V= K[\widebar E]$. Now define the ordering on $\power(\kappa)$ given by:
%$$ x \triangleleft y \mbox{ iff } \exists M( M \mbox{ an $\widebar E$-mouse with }\rho^\omega_M=\kappa \wedge x<_M y)$$
%where $<_M$ is the usual order of the structure of $M$. 
That this is a good $\Sigma_1(\kappa, e)$-wellorder follows directly from:\\

\nod{\em Claim: If $M,N$ are two mice satisfying the above for $x,y\sset \kappa$ then $M=N$.
}
%with $\rho^\omega_M=\kappa = \rho^\omega_N$ then $x<_Ny \equi x<_M y$.
%}
%The Claim follows from the:
%\nod{\em Subclaim: Any two such $M,N$ one is an initial segment of the other
%} since then there can be no disagreement about the ordering of $x$ and $y$.\\

\pf of Claim: by standard comparison considerations, which we shall give in any case. Let $M=M_0$ and $N= N_0$ be two such mice; let them be compared to $M_\theta$ and $N_\theta$. We want that the comparison is trivial, \ie $M=N$.  Suppose for a contradiction that $\nu_0$ is the point of least difference between $E^{M_0}$ and $E^{N_0}$. As they both satisfy  $\psi(E,\lambda)$ there are no measurable cardinals in the interval $(\lambda^{+},\kappa)$ %of $\bar K\dfs K[\widebar E]$, interval 
on either of the $E^{M_0}$, $E^{N_0}$ sequences.  Suppose first that $\psi(E,\lambda)$ asserts only that the measurables are bounded by $\lambda$, that is the measurables (and their measures) in  $M,N$ are just those in $e$. Thus were $\nu_0<\kappa$ we should have a truncation on one side, wlog, the $N$-side to create a full measure to form an ultrapower. As there can be no truncation on the $M$-side by a consequence of the Dodd-Jensen Lemma, the comparison must run for at least $\theta \geq\kappa$ stages finally iterating some measure of order zero in some $N_\iota$ up to $\kappa$. As there are no measures in this interval on the $M$-side to take ultrapowers with, then there has been no movement on the $M$-side: $M_0=M_\kappa$. 
However the iteration of the initial truncate $N^\ast_0$ of $N_0$ to $N_\kappa$ is from some stage before $\kappa$ onwards,  a simple iteration (after perhaps finitely many further truncations) that can be defined inside the $KP$ model $N_0$. We may conclude that $N_\kappa\in N_0$. $N_\kappa$ is of the form $(J^{E^{N_\kappa}}_\alpha, F^{N_\kappa})$ for some filter $F^{N_\kappa}$. Now $M=M_\kappa$ is a simple $KP$ model, with $\kappa$ as its largest cardinal. Hence it is a proper initial segment of the $ZF^-$-model $J^{E^{N_\kappa}}_\alpha$ and thus is an element of $N$. But this contradicts the assumption on the $\in$-minimality of $N$.
Consequently any non-trivial comparison must start by using some $\nu_0 >\kappa$ indexing some filter with critical point $\geq \kappa$. However this is also a contradiction since both $\rho^\omega_M=\kappa =\rho^\omega_N$, our %other minimality 
conditions insure that if $M\neq N$ then we see by comparison that the code of one as a subset of $\kappa$ is a member of the other. But that also contradicts the minimality conditions on the appearance of $x,y$ in the two hierarchies above $\kappa$.  We conclude that $M=N$.

In the case that in $E^{K}$ that $\lambda$ is strong up to $\kappa$ then let $M$ be some initial admissible segment $K$ satisfying the requirements. 
%Let $F_{0}$ be this extender on $E^{M}$. 
Suppose $N$ is another mouse satisfying them with $\lambda$ strong up to $\kappa$
% as witnessed by an extender $F_{1}$.
%But $\tau \neq \tau^{N}$ would engender that $0^{pistol}$ exists which is absurd. 
But   the extenders on the $E^{N}$ sequence must agree with those on the $E^{K}=E^{M}$ sequence below $\kappa$. Otherwise in the comparison of $M$ with $N$ if $\nu_{0} $ is the least index used, this must be because both $F_{0}\dfs E^{M}_{\nu_{0}} $ and $F_{1}\dfs E^{N}_{\nu_{0}}$ are both non-empty. But we are in $K$ and $\power(\lambda)\in M\cap N$. Thus both $F_{0},F_{1}$ are $\omega$-complete \cite{Z02} Lemma 8.2.12); this guarantees that $\la J^{E^{K}}_{\nu_{0}},\in, E^{K},F_{0},F_{1}\ra $ is a bicephalus. 
%(See the discussion in \cite{Z} Ch. 8 and Lemma 
And thus $F_{0}=F_{1}$ ({\em op.cit.} Lemma 8.2.9).
% is on the $E^{K}$-sequence.
Thus if any comparison is to be done it must involve an index $\nu_{0}>\kappa$ indexing an extender with critical point $>\kappa$ (by $\neg 0^{pistol}$). But just as before this contradicts our minimality conditions on $M,N$ and we conclude that $M=N$.
%there is on the $E^{N}$ sequence an $\omega$-complete extender, strong up to $\kappa$, which is different from its counterpart on the $E^{K}$ sequence.  But this is impossible. So again we conclude that $M=N$.
%For, if in the comparison there was a truncation on one side to create a full measure, this ensures that $0^{pistol}$ would exist.%by our requirements 
%By our assumption $(\ast)$ $\tmop{crit}(E^{N_{\nu_0}}) \geg\kappa$ (as $\kappa$ is not overlapped). Thus if $\nu_0$ is defined it must be that of a truncation of $N=N_0$ to a $ZF^-$ model $J^{E^N}_\alpha$, with a final topmost measure $F$ with critical point $\kappa$ and
%with $E^N\rest\alpha = E^M\rest\alpha$. Then, again, $J^{E^N}_\alpha\models ZF^-$ whereas by assumption $J^{E^M}_\alpha\models KP$ but by the $\in$-minimality clause, no stronger theory. So this is absurd. We conclude the coiteration is trivial and $M=N$.
% $\power(\kappa)^M = \power(\kappa)^{N_0^\ast}$.
% $N_\kappa$. 
%But $M_0 <^* N_\kappa$ in the mouse order, as $ 
%\rho^\omega_M=\kappa > \rho^\omega_N$. But then $A_M$ the canonical code for $M$ is an element of $N_\kappa$ and hence by the internal definability of the iteration so far, is an element of $N$.
\qed (Claim and Lemma)\\

\begin{corollary}($V=K$)
Let $\kappa$ satisfy the assumption $(*) $ of the last lemma. Then 
$$ \kappa \mbox{ has the \siclub \, property } \Equi \kappa \mbox{ is stably measurable}\Equi \sigma ( \kappa ) =u_{2} ( \kappa ) .$$ 
\end{corollary}
\pf We just repeat as before that $\power(\kappa)$ having a good 
{\boldmath$\Sigma_1$}$(\kappa)$-wellorder together with the  \siclub \, property   implies that $(\widetilde M,\in , F_{\kappa})$ witnesses stable measurability. The right-to-left direction of the first equivalence is now trivial. The second equivalence is Lemma \ref{u=sigma} as before.\mb \qed\\

As we saw at Lemma \ref{L1.27}, without an assumption the first equivalence can fail, for example $\kappa$ a regular limit of measurables, which is not a Mahlo cardinal.
%\todoi{Under the assumption we should have $u_2= \sigma$, using the same arguments as before.}

\begin{lemma} Assume $\neg 0^{pistol}$.
$\sigma ( \kappa ) =u_{2} ( \kappa )   \Imp$ $\sigma ( \kappa )=\sigma ( \kappa )^{K}
  =u_{2} ( \kappa )^{K} $.
  % If additionally $\neg 0^{\dagger}$ \ then $\sigma ( \kappa ) = \sigma ( \kappa )^{K}$.\\

\end{lemma}
\pf The assumption implies that bounded subsets of $\kappa$ are closed under $\sharp$'s. By $\neg 0^{pistol}$ and absoluteness arguments $u_{2}= u_{2}^{K}$. \qed

\subsection{When $K=K^{{JS}}$} 

 Let $V=K^{{JS}}$ be the Jensen-Steel core model built assuming there is no inner model of a Woodin cardinal. Then the comparison argument in Lemma \ref{strongwo} goes through with the same effect, for a $\kappa$ which is not a limit of $K$-measurables.
 %or is such that there are many strongs up to $\kappa$. However this says little: there may be only some small $\kappa$ satisfying this.
 %  $(*)$, if $K$ has strong cardinals.

%\todoi{This needs changing/deleting}

\section{Two applications}
\subsection{$\Sigma_{1}$-stable measurability}
There are two further recent theorems that could benefit from the weakening
of an assumption from iterability to stable measurability. They are proven in
{\cite{LSS18}} as Theorems 1.9 and 1.8 respectively with the assumption of
($\omega_{1}$-)iterability, which we now weaken to stable measurability by
adapting their argument. But the proofs are now shorter.

\begin{theorem}\label{2.1}
  Assume $\kappa$ is stably measurable. Then the following are equivalent for
  $X \sset \re$:
  
  (i) $X$ is $\Sigma_{1} ( \kappa )$-definable; \ \ (ii) $X$ is
  $\Sigma^{1}_{3}$ definable.
\end{theorem}

{\bf Proof:}
  (ii) $\Imp$ (i) is unaltered as in {\cite{LSS18}}. \ (i) $\Imp$ (ii): Let
  $\pa{\widetilde M, \in ,F_{\kappa}}$ witnesses stable measurability. Exactly as in
  {\cite{LSS18}} define the $\Sigma^{1}_{3}$ set:
  $$Y= \left\{ y \in \re \mid \ex \mbox{ countable, iterable } \pa{N, \in ,U} \models \right.\mbox{``}U \mbox{ is a
  normal ultrafilter on } \bar{\kappa}\wedge \varphi ( \bar{\kappa} ,y)\mbox{''}  \}$$
  
  {\nod}where $\varphi ( \kappa ,v_{1} ) \in \Sigma_{1}$ defines $X$.  Then $Y \supseteq X$ since for any $y \in X$ we can
  take a countable elementary substructure $\pa{N_{0} , \in ,U_{0}} \prec
  \pa{\widetilde M, \in ,F} \models \varphi ( \kappa ,y )$. Then we have a witness to put
  $y$ into $Y$. Conversely  any witness $\pa{N_{0} , \in ,U_{0}} \models
  \varphi ( \bar{\kappa} ,x )$ that $x \in Y$, iterates to a structure \
  $\pa{N_{\kappa} , \in ,U_{\kappa}} \models \varphi ( \kappa ,x )$, with
  $N_{\kappa} \in \wq$. But $\wq = \widetilde M$ by Lemmata \ref{Lemma1.10} and  \ref{1.11}. But then by $\Sigma_{1}$-upwards absoluteness $  \varphi ( \bar{\kappa} ,x )$ holds in ${\widetilde M}$ and in $V$.
\qed\\

For completeness we repeat the following immediate, but nice, corollary 6.3 from \cite{LSS18} with this improved hypothesis.
\begin{corollary} Suppose $\kappa$ is stably measurable. Then if there is a $\Sigma_{1}(\kappa)$ wellordering of $\re$ then there is such which is $\Sigma^{1}_{3}$.
\end{corollary}

In $K$ we get a form of equivalence in Theorem \ref{2.1}.

\begin{theorem}\label{2.3} Assume $V=K^{DJ}$.
% and that $H(\kappa)$ is closed under sharps.  
Let $\Phi(\kappa)$ be the following sentence:
$$ \Phi(\kappa) : \forall X \sset \re\all r \in \re[
  X\mbox{ is }\Sigma_{1} ( \kappa, r )\mbox{-definable} \equi X\in \Sigma^{1}_{3}(r)].$$
  Then we have:
 $$\kappa\mbox{ is stably measurable }\equi  \Phi(\kappa) \mbox{ is preserved by small forcings of size } < \kappa.$$

\end{theorem}

 \pf By Theorem \ref{2.1} $\kappa$ is stably measurable implies $\Phi(\kappa)$, and stable measurability is preserved by small forcing.  This proves $(\rightarrow)$.
 %we only need to show $(\leftarrow)$. We have the following facts.

%{\cor [(1)  The first $\omega$ iteration points of $Q(\kappa)$ are the uniform indiscernibles $\pa{u_{n}\mid 1<n<\omega}$. ]}
First just note that if $H(\kappa)$ is not closed under sharps (which implies that $\kappa$ is not stably measurable) then the right hand side fails: let $a\sset \gamma < \kappa$ have no sharp; then $V=L[a]$ as we must have $K\sset L[a]$ for otherwise $a^\sharp$ would exist. Let $\mathbbm{P} = Col(\omega, \gamma)$; then $V[G]\models V= L[r]$, where $r$ is a real coding $G$ and $a$. But now any analytical (in $r$) set whatsoever is definable over $L_{\omega_1}[r]$ and thus is $\Sigma_1(L_\kappa[r],r)$ and then $\Sigma_1(\kappa,r)$. Consequently the right hand side fails. 

So now assume that $H(\kappa)$ is  closed under sharps.

\nod(1) Any $ X\in \Sigma^{1}_{3}(r)$ is $\Sigma_{1}^{Q(\kappa)}(r)$. This follows from the fact that there is a Martin-Solovay tree for $\Pi^1_2$ is $\Delta_{1}$-definable over $Q(\kappa)$ (cf. \cite{We86}, \cite{W11} Sect 1).

%\todoi{More detail here perhaps}
\nod (2) %$\widetilde M$ is a transitive set such that $\{\kappa\} \cup H_{\kappa}\sset 
$\widetilde M\prec_{\Sigma_{1}}H(\kappa^{+})$, by Lemma \ref{1.28}
% is any set $\Sigma_{1}$-elementary in $H(\kappa^{+})$, 
and then by definition $\sigma\leq On\cap \widetilde M$.

 Suppose $\kappa$ is not stably measurable. Then $Q(\kappa)$ cannot witness stable measurability and moreover:
  
\nod (3) $\theta(\kappa)\leq u_{2} < \sigma$. \\
 \pf The first inequality is Lemma \ref{utheta}, and the second  is by Theorem \ref{ordinalsinK}. \qed (3)\\
 %$u_{2}=\theta(\kappa)$ and 

 But then: 
 
\nod (4) $Q(\kappa)\in \widetilde M$. \\
 \pf We have that $\theta(\kappa)= On \cap Q(\kappa) < \sigma\leq On \cap \widetilde M$. But then for some $z\in H(\kappa)$, $\theta(\kappa) \in \Sigma_1^{\widetilde M}(\kappa, z)$. But then also $J_{\theta(\kappa)}^{F_\kappa}$ is also $\Sigma_1^{H(\kappa^+)}(\kappa, z)$, and  so is in ${\widetilde M}$. \qed (4)\\
 
 Let $G$ be $\mathbbm{P}$-generic over $V$ for some $\mathbbm{P}\in H(\kappa)$ which collapses $TC(\{z\})$ to be countable. Then  as
 $\widetilde M\prec_{\Sigma_{1}}H(\kappa^{+})$, we have in $V[G]$
 $\widetilde M[G]\prec_{\Sigma_{1}}H(\kappa^{+})[G] = (H(\kappa^{+}))^{V[G]}$. 
 Let $r\in \re^{V[G]}$ code $z$. Then $Q(\kappa)$, which is not altered in the passage to $V[G]$, is  
 in $\Sigma_{1}^{\widetilde M}(\kappa, r)$.
 
 %$\widetilde M$ is a $\Sigma_{1}$-elementary substructure of $H(\kappa^{+})$ and every element of $\widetilde M$ is $\Sigma_{1}(\kappa)$ definable (in $\widetilde M$, in $H(\kappa^{+})$, or indeed $V$). 
 %\todoi{NO: and every element of $\widetilde M$ is $\Sigma_{1}(\kappa)$ definable in parameters!}But then $\theta(\kappa)$ is $\Sigma_{1}(\kappa)$-definable.
 % as the largest $\theta < \xi$ so that $J_{\theta}^{F_{\kappa}}\models  $``$F_{\kappa}$ is a normal ultrafilter on $\kappa$'' but $J_{\theta+1}^{F_{\kappa}}$ does not. As $\xi$ is $\Sigma_{1}(\kappa)$ definable, so is $\theta(\kappa)$ and 
% Thus so is $ J_{\theta(\kappa)}^{F_{\kappa}}$ and this places it in $ \widetilde M$. But the latter is $Q(\kappa)$. \qed (4)

 Consequently if $X\sset \re$ is a  universal $\Pi^1_3$ set, then  $X\in \Pi_{1}^{Q(\kappa)}$ but would then be  $\Sigma_1^{\widetilde M}(\kappa, r)$; but such an $X$ is not $\Sigma^1_3(s)$ for any $s\in \re$.  So this provides a counterexample to the preservation of $\Phi(\kappa)$ under small forcing.
  % but is nevertheless in $\widetilde M$ and so $\Sigma_{1}(\kappa)$ definable. $X$ is thus a counterexample to the equivalence on the right hand side of the theorem's statement. \\ \mb\hfill
 %{\em Case 1 $Q(\kappa)$ is admissible.} 
  \qed \\
  
  Within $K$ we can replace the stable measurability by any of its equivalents from Theorem \ref{ordinalsinK} of course. Outside of $K$ even assuming sufficient sharps for $\Sigma^{1}_{3}$-absoluteness  we can only show by similar methods results such as the following:
  
  \begin{lemma}\label{3.4} $(\neg 0^{dagger} \wedge \forall a\in \power_{<\kappa} (\kappa)(a^{\#}$ exists$))$. Assume there is a good {\boldmath$\Sigma_1$}-wellorder of $\power(\kappa)$. Then:
  $$ u_{2}(\kappa)<\sigma(\kappa) \Imp  \Phi(\kappa)
 % X\in \Sigma_{1} (r) \equi X\in \Sigma^{1}_{3}(r) 
  \mbox{ fails in a small generic extension}.$$
  
  \end{lemma}
  \pf  Use that if $\widetilde M$ is a $\Sigma_{1}$-substructure, that  $\theta(\kappa) \leq
  (u_2)^{K} = u_{2}$ by the correctness of the calculation of $u_{2}$ inside $K$ due to the assumed absoluteness from $\neg 0^{dagger}$, and  thus is $\Sigma_{1}^{\widetilde M}(\kappa, z)$ definable from some $z\in \power_{<\kappa} (\kappa)$, and thus also $Q(\kappa)\in \widetilde M$ as above. But then the first $\omega$ iterates of $Q$ are in $\widetilde M$ and this is enough to define the   Martin-Solovay tree of $K$ on these uniform indiscernibles as an element of $\widetilde M$.  (The assumptions of the lemma again ensure the correctness of this tree in $V$.)
  But now we get as before $\Pi^{1}_{3}$ sets of reals as $\Sigma_{1}(\kappa, r)$ where $r$ is a real in a small generic extension coding $z$.
  \qed \\
  
  But we don't have a converse to this.
  
%\todoi{Under the special assumption, if $\kappa$ is not s.m. then we'll get $\wb \cap On = u_2 < \sigma = On \cap \widetilde M$. (But do we get $\wb \in \widetilde M$?) Is this enough to get $T_{MS} \in \widetilde M$?}

\begin{theorem}
Assume $V=K^{strong}$. Let $\kappa$ not be a  limit of measurable cardinals. Then the conclusion of the last Theorem \ref{2.3}  holds.

\end{theorem}
\pf The direction $(\rightarrow)$ is as before, again we seek to prove $(\leftarrow)$. Instead of using the Dodd-Jensen $Q(\kappa)$ we use the generalised $\wq(\kappa)$. 
 If $\wq(\kappa)$ is in $\widetilde M$ we'll reason as before that if $\widetilde M$ fails to witness stable measurability, that analytical sets are definable over $\wq(\kappa)$ because again a Martin-Solovay tree is so definable.
 We again then have a counterexample to the right handside. 
 
 The case that $H(\kappa)$ is not closed under sharps is a small variant: let $a\sset \gamma < \kappa$ have no sharp; let $a'$ code both $a$ and $K\rest \gamma'$ where $\gamma' < \kappa$ is least with $a\in K\rest\gamma'$. Then $V=L[a']$. Let $\mathbbm{P} = Col(\omega, \gamma')$; then $V[G]\models V= L[r]$, where $r$ is a real coding $G$ and $a'$. We can finish as before, 
 %But now any analytical (in $r$) set whatsoever is definable over $L_{\omega_1}[r]$ and thus is $\Sigma_1(L_\kappa[r],r)$ and then $\Sigma_1(\kappa,r)$. Consequently the right hand side fails. 

 We assume then $H(\kappa)$ is closed under sharps; we are done if we can show $\wq(\kappa) \in \widetilde M$. Note that by Lemma \ref{strongwo}   we have a good $\Sigma_1(\kappa, e)$ wellorder of $\power(\kappa)$ and hence 
$\widetilde M\prec_{\Sigma_{1}}H(\kappa^{+})$. (Recall that $e$ was the initial segment of the $E^{K}$ extender sequence $E^K\rest \lambda^{+}$ for some $\lambda < \kappa$ which bounds the measurable cardinals.)
By the assumed failure of stable measurability at $\kappa$ we must have $\widetilde M \neq \wq(\kappa)$ as otherwise $(\widetilde M, F^\kappa)$ would be a witness to this. Let $A\in M_0$ be such that $A\in \widetilde M\back \wq(\kappa)$. Without loss of generality we may assume $A\rest \lambda^{+}$ codes $e = E^K\rest \lambda^{+}$.

Firstly suppose that $\neg A^\sharp$. Then covering lemma arguments show that $K^{A}\dfs (K)^{L[A]}$ is a universal weasel, and as we are below $0^{pistol}$ it is a simple iterate of the true $K$ - that is without truncations in the comparison. However $A$ codes the initial segment of $K$ given by $E^{K}\rest \lambda^{+}$ and thus $E^{K^{A}}\rest  \lambda^{+} =E^{K}\rest \lambda^{+}$.
Consequently no comparison index is used below $\kappa$. 
Consequently we have that $K^{A}_\kappa = K_{\kappa}= L_{\kappa}[A]= H(\kappa)$.
%As $V=K$ we must have that $L[A] = K$. In particular that $L_\kappa[A]= K|\kappa = H(\kappa)$. 
But $L_\kappa[A]\in \widetilde M$. But then $\wq(\kappa)$ is definable within the admissible set $\widetilde M$ from $H(\kappa)$ and we've achieved our goal.

Thus suppose $A^\sharp$ exists. If $L_\kappa[A]= K_\kappa =H(\kappa)$, then we could reason as we just have done that  $\wq(\kappa)$ is definable within $\widetilde M$. So there is some $<^\ast$-least sound mouse $P$ with 
$A\in P$ and $\rho^\omega_P=\kappa$. By the elementarity of $\widetilde M$ in $H(\kappa^+)$ we have that $P\in \widetilde M$ as it is $\Sigma_1$ definable from $A$.
%$P\in H(\kappa)\back L_\kappa[A]$. 
Then in comparison of $P=P_0$ with $R_0 \dfs K_ \kappa$ we cannot have that $R_0$ is truncated below $\kappa$ and some $R_0^\ast$ is iterated past $P$, as in that case $A$ is an element of an iterate of the $\kappa$'th iterate of (some final truncate of) $R_0^\ast$, and the latter along with $A$ would be in $\wq(\kappa)$. 
%But equally, no $<\kappa$-sized initial segment of $R_{0}$, say $K\rest \tau$  for a $\tau <\kappa$, can iterate past $P$ either.
So then, as $K$ has no full measures in the interval $(\lambda,\kappa]$, the coiteration is trivial below $\kappa$, indeed altogether trivial, and $H(\kappa) = K\rest \kappa \in P\sset \widetilde M$, and we may finish as before.
\qed\\

Putting together the above we have: 
\begin{corollary} Assume $V=K^{DJ}$ (or $V=K^{strong}$). % and is closed under sharps. 
Then $\exists \kappa \Phi(\kappa)$ is (set)-generically absolute if and only if there are arbitrarily large stably measurable cardinals in $K$.
\end{corollary}

 % \pf  We just consider the second case of  $V=K= K^{strong}$ as the other is simpler.  Suppose there are arbitrarily large stably measurable cardinals. Then no cardinal is strong in $K$: If $\lambda$ were strong and $\kappa > \lambda$ stably measurable, then every bounded subset of $\kappa$ has a sharp, and this would imply $0^{pistol}$ exists in $K^{strong}$ which is absurd. Then there are arbitrarily large stably measurable cardinals. (Because there  is  a proper class of non-measurable Ramseys which are themselves not limits of measurables and each of which is {\em a fortiori} stably measurable.)
 % \qed\\
 
 As in Lemma 3.4 we can prove the following with these methods.
 
 \begin{corollary} Assume $\neg 0^{pistol}
  \wedge \forall a\in \power_{<\kappa} (\kappa)(a^{\#}$ exists$))$. Assume there is a good {\boldmath$\Sigma_1$}-wellorder of $\power(\kappa)$. Then:
  $$ u_{2}(\kappa)<\sigma(\kappa) \Imp  \Phi(\kappa)
 % X\in \Sigma_{1} (r) \equi X\in \Sigma^{1}_{3}(r) 
  \mbox{ fails in a small generic extension}.$$

 % Then $\exists \kappa \Phi(\kappa)$ is (set)-generically absolute if and only if there are arbitrarily large stably measurable cardinals in $K^{strong}$.
 \end{corollary}

 The following is a strengthening of \cite{LSS18} Theorem 1.8 where the assumption is that $\kappa$ is iterable; it is based on their template but now follows easily from the analysis above.
 
\begin{theorem}
  Assume $\kappa$ is stably measurable. Assume $X \sset \power ( \kappa )$
  separates $F_{\kappa}$ \ from $\tmop{NS}_{\kappa}$, then $X $ is not
  $\Delta^{H ( \kappa^{+} )}_{1}$.
\end{theorem}

{\bf Proof:}
  Let $\kappa$ be stably measurable as witnessed by $\pa{\widetilde M, \in ,F_{\kappa}}$ as usual. %By Corollary \ref{Cor7} we may assume $F=F_{\kappa} \cap M$. 
  For a contradiction let $\varphi (
  v_{0} ,v_{1} )  $ and $\psi ( v_{0} ,v_{1} )  $ be  $\Sigma_{1}$ and define
  some $X \supseteq F_{\kappa}$ and its complement in $\power ( \kappa )$, but
  with $X \cap \tmop{NS}_{\kappa} = \emp$. Then $F_{\kappa} \cap \widetilde M$ is
  $\Delta^{\widetilde M}_{1}$ and the statement that it is an ultrafilter is
  $\Pi_{1}^{\widetilde M}$. As $\widetilde M \prec_{\Sigma_{1}} H ( \kappa^{+} )$, we thus have an
  $\widetilde{F } \supseteq F$, which is a definable $H ( \kappa^{+}
  )$-ultrafilter. But this is absurd, as {\cite{LSS18}} \ says, as then
  $\widetilde{F}  $ is definable over $H ( \kappa^{+} )^{\tmop{Ult} ( H (
  \kappa^{+} , \tilde{F} ) )}$. 
\qed

\section{When $\sigma > u_{2}$ and canonical models}

The following definition can be given a first order formulation as  a scheme.
\begin{definition}
Let $\varphi(v_{0})$ be a formula of the language of set theory with the single free variable $v_{0}$. Let $M$ be an inner model of ZFC (thought of as a transitive proper class of sets defined by some class term). We say that $M$ is {\em canonically defined by} $\varphi(\xi)$ (for some parameter $\xi \in On$), if $\varphi(\xi)^{M}$ but for no other inner model $M'$ do we have $\varphi(\xi)^{M'}$. 
\end{definition}
Clearly then $L$ is such (``$V=L$'') but also $L[\mu]$ (``$V=L[\mu]$ where $\mu$ is a normal measure on $\kappa$'' - using the ordinal parameter $\kappa$. ``$V=K$'' by itself does not canonically define any inner model, but $L[0^{\#}]$ or the least inner model where all sets have $\#$'s, $L^{\#}$, are canonical in this sense. Hence Carl and Schlicht ask: what is the least $L[E]$-model which is not canonical? Clearly if an inner model thinks that it is not canonically definable, then it is a model of an inner model reflection principle (see Def. \ref{IMR} below). Then \cite{BCFHR18} ask for upper bounds to the existence of a model of inner model reflection, thus essentially the same question.

We identify this model, as an inner model, and it turns out to be an inner model of the full Dodd-Jensen core model below a measurable cardinal.
It is a model which is intermediate in consistency strength between admissible measurability and stable measurability. In this model no $Q$-structure $Q ( \kappa ) =
\pa{J^{F_{\kappa}}_{\theta ( \kappa )} , \in ,F_{\kappa}}$ witnesses stable measurability, but such can be admissible, and moreover can be first order reflecting.

\begin{definition} A transitive admissible set $\mathcal{A}$ is {\em first order } (or $\Pi^{0}_{\omega}$) {\em reflecting} if for any formula $\varphi(\vec p)$ with parameters $\vec p \in \mathcal{A}$ such that $(\varphi(\vec p))^{\cala}$  there is a transitive $u\in \cala$ so that
$(\varphi(\vec p))^{u}$. 
\end{definition}
We shall adopt a version of this appropriate for $Q$-structures: for $u$ we just take a proper initial segment of $Q$. (Note that ``$V=L[F]$'' is in any case $\Pi_2$ so this is not a restriction.) 
\begin{definition}
$Q(\kappa)$ is $\Pi^{0}_{\omega}$ {\em reflecting} if for any $\varphi(\vec p)$ with parameters $\vec p \in H(\kappa)$ with $Q(\kappa)\models \varphi(\vec p)$ there is $\tau$ with $\kappa \leq
\tau < \theta(\kappa)$ with $J^{F_{\kappa}}_{\tau}\models \varphi(\vec p)$.
\end{definition}

We shall tie this up with a version of: 

\begin{definition}[{\em Inner model reflection}]\label{IMR}(i) An inner model $M$ is {\em reflecting for}  $\varphi(p)$, for $p\in
M$, $\varphi \in \call_{\dot \in, \dot =}$ when, if it is a model of $\varphi(p)$, then there is a proper inner model $M'\subset M$ which is a model of $\varphi(p)$.\\
(ii) An inner model is {\em first order reflecting} if it is first order reflecting for all $\varphi(p)$. It is {\em $\Pi_{n}$-reflecting} when it is so for all $p$ and all $\psi(v_{0})\in \Pi_{n}$.
\end{definition}
%(We use finite sequences of ordinals in $\vec\xi$ rather than set parameters, so that we do not have to fuss about whether parameter sets are in the proper inner models \etc) 
Clearly a model which is first order reflecting cannot be canonical in the sense above.

Given a mouse $N$ (in the modern sense) in $K^{DJ}$ this generates an inner model $K^{N}$ (which is of the form $L[E^{K^{N}}]$ for some predicate $E^{K^{N}}$).  Let $C_{N}= \la \kappa_\alpha \mid \alpha \in On \ra$ be the cub class of the iteration points of $N$ as we iterate by its top active measure. It then generates the inner model $K^{N}=\bigcup _{\alpha\in On}H_{\kappa_{\alpha}}^{N_{\alpha}}$. As above we can let $Q^{N}(\gamma)$ be the union of all the Dodd-Jensen mice in $H_{\gamma}^{K^{N}}$ iterated to comparability at $\gamma$.

\begin{theorem} Let $N$ be a mouse that generates an inner model $K^{N}$ which is $\Pi_{n}$-reflecting.  Then for any $\kappa \in C_{N}$, $Q^{N}(\kappa)$ is $\Pi^{0}_{n}$-reflecting.
Conversely if $N$ generates the model $K^{N}$ so that  for some $\kappa \in C_{N}$, $Q^{N}(\kappa)$ is $\Pi^{0}_{n}$-reflecting, then $K^{N}$ is $\Pi_{n}$-reflecting.
\end{theorem}

\pf %We assume $n>1$. 
%Suppose $K^{N}$ is $\Pi_{n}$-reflecting. 
Recall that for any $\gamma\in \tmop{Card}^{K^{N}}$, $K^{N}_{\gamma}= (H({\gamma}))^{K^{N}}$. As the elements of $C_{N}= \{\kappa_{\alpha}\mid \alpha\in On\}$ are indiscernible for $K^{N}$ we shall have that for any $\alpha$ so that $\kappa_{\alpha} > \max \tmop{rk}_{K^{n}}(\vec p)$:
$$(1) \quad\quad\pa{{K^{N}},\in}\models \varphi(\vec p)\, \equi \, \pa{{K^{N}_{\kappa_{\alpha}}},\in}\models\varphi(\vec p)\, \equi \,\pa {Q^{N}(\kappa_{\alpha}),\in}\models \mbox{ ``}\pa{H(\kappa_\alpha),\in}\models  \varphi(\vec p)\mbox{''}.$$
%noting  

For a $Q(\kappa)$-mouse the first  projectum $\rho_{Q}^{1}$ is $\kappa$  (indeed all projecta are). By the fine structure for such mice, we have that any $\Pi_{n}^{\la Q(\kappa),\in \ra}$ relation $R\sset H(\kappa)^{Q(\kappa)}$ is 
$\Pi_{n}$ over $\la H(\kappa)^{Q(\kappa)},\in\ra$.
(Officially because we use that $J_{\rho^{1}_{Q}}^{{A^{1}_{Q}}} = H(\kappa)^{Q(\kappa)}$, where $A^{1}_{Q}$ is the first mastercode of $Q=Q(\kappa)$. We shall write $Q^{N}(\kappa)$ for $Q^{K^N}(\kappa)$.) Using this with $\kappa = \kappa_{\alpha}$ in the equivalences at (1), together with  $H(\kappa_\alpha)^{K^{N}} = K^{N}_{\kappa_{\alpha}}= H(\kappa_\alpha)^{Q^{N}(\kappa_\alpha)}$ we have the equivalence of the right hand statement with $ \pa {Q^{N}(\kappa_{\alpha}),\in}\models \mbox{ ``} \varphi(\vec p)\mbox{''}$.

%$\Pi_{n-1}$ over $\la H(\kappa)^{Q(\kappa)},\in, A^{1}_{Q(\kappa)}\ra$, where $A^{1} _{Q(\kappa)} = \{ \pa{ i, x}\mid Q\models \varphi_{i}(x)\}$.  (Here $\pa{ \varphi _{i}\mid i \in \omega}$ is a recursive enumeration of the $\Sigma_{1}$ formulae with a single free variable: $\varphi(v_{0})$. Standard parameters for such $Q$ are all empty and so are uninvolved in these definitions.)

Now suppose $K^{N}$ is $\Pi_{n}$-reflecting, there is an inner model $K'\subset K^{N}$ in which $\varphi(\vec p) + \mbox{``$V=K$''}$ holds (at least if $n\geq 2$; if $n=1$ it reflects to $L[\vec p]$). But any such model $K'$ is actually some $K^{M}$ for an $M<^{\ast}N$ (and thus $M\in K^{N}$). Now choose $\alpha$ sufficiently large so that it is greater than $|M|^{K^{N}}$ and is also in $C_{M}$. As $M$ is missing from $K^{M}$ it is easy to see that $Q^{M}(\kappa_{\alpha})$ is a proper initial segment of $Q^{N}(\kappa_{\alpha})$. However the sequence of equivalences in (1) holds with $M$ replacing $N$ throughout.

Now for the converse suppose $K^{N}\models \psi(\vec p)$, and then \via (1) above, we have $Q^{N}(\kappa)\models \psi(\vec p)$ for a $\psi\in \Pi_{n}$, $\kappa\in C_{N}$, and $Q^{N}(\kappa)$, $\Pi^{0}_{n}$-reflecting.
We note that $Q^{N}(\kappa) = \bigcup \{M_{{\kappa}} \mid {M \in  Q^{N}(\kappa)},\,\tmop{On}\, \cap\, M < \kappa, \,M \mbox{ a $DJ$-mouse}\}$; the latter is described by a $\Pi^{0}_{2}$ formula, which may assume then is a conjunct of the formula $\psi$. (Another way of putting this is to say that $\theta^{{N}}({\kappa})$ is a multiple of $\kappa$.)
By (1) again we have  this is equivalent to  
$K^{N}_{\kappa}\models \psi(\vec p)$. As $Q^{N}(\kappa)$ is $\Pi_{n}$-reflecting, there is some ${\kappa}<\tau < \tmop{On}\,\cap\, Q^{N} (\kappa)$ with $J_{\tau}^{F_{\kappa}}\models \psi(\vec p)$ and $J_{\tau}^{F_{\kappa}} = \bigcup \{M_{{\kappa}} \mid {M \in J_{\tau}^{F_{\kappa}} }, \tmop{On}\, \cap\, M < \kappa, M \mbox{ a $DJ$-mouse}\}$. Now with this property of $\tau$, this ensures that the $<^{\ast}$-least mouse $M\notin H^{J_{\tau}^{F_{\kappa}}}_{\kappa}$ with $\tmop{crit}(M)=\kappa$ generates a proper inner model of $K^{N}, K^{M}$, with $Q^{M}_{\kappa}=J_{\tau}^{F_{\kappa}}$ and now, applying (1) once more, $(\psi(\vec p))^{K^{M}} $.
%We claim there is an effective procedure to produce a $\psi'$ so that $Q^{N}(\kappa_{\alpha})\models \psi(\vec \xi)$ iff $H(\kappa_{\alpha})^{K^{N}}\models \psi'(\vec\xi)$. We note that $H(\kappa_{\alpha})^{K^{N}} \sset Q^{N}(\kappa_{\alpha})$ and the former is a definable sub-class of the latter (without being an element). As an example suppose $\psi$ is $\ex u \all v \chi (u,v,\vec\xi)$. [We may assume that it is of the form $\ex \tau \all \tau' > \tau J_{\tau'}^{F_{\kappa}}\models  \chi (u,v,\vec\xi)$]
%But any $u \in  Q^{N}(\kappa_{\alpha})$ is in fact an element of some $M^{0}_{\kappa_{\alpha}}$ for an $M^{0}\in H(\kappa_{\alpha})^{K^{N}}$. And to express ``$\all v$'' over $Q^{N}(\kappa_{\alpha})$ this will be equivalent to ``$\all M^{1}\in H(\kappa_{\alpha})^{K^{N}} \all v \in M^{1}_{\kappa_{\alpha}}$.''
\mb\hfill \qed \\

Perhaps unsurprisingly, there is a strict hierarchy under $\subset$ of $\Pi_{n}$-reflecting inner models in $K^{DJ}$ for increasing $n$.

\begin{corollary} For $n > 1$, if a mouse $N$ generates an inner model $K^{N}$ which is $\Pi_{n}$-reflecting, then for $\kappa\in C_{N}$, we have that $Q^{N}(\kappa)$ is a $\Pi^{0}_{n}$-reflecting admissible set. Furthermore for such $\kappa$ there is $M <^{\ast} N$ and a $U^{M_{\kappa}}$ measure one set of $\xi < \kappa$ such that $Q^{M}(\xi)$ is $\Pi^{0}_{n-1}$-reflecting. Hence $K^{M}$ is a proper inner model of $K^{N}$ which is $\Pi_{n-1}$-reflecting.

\end{corollary}
\pf The first sentence is just a restatement of part of the proof above.
The statement ``$Q({\kappa})$ is $\Pi^{0}_{n-1}$ reflecting'' is itself a $\Pi^{0}_{n}$ statement over $Q(\kappa)$:
$$ \forall \varphi \in \Pi_{n-1}\all x \in H_{\kappa}[\varphi(x)\Imp \exists \tau J_{\tau}^{F_{\kappa}}\models \varphi(x)].$$
%is $\Pi_{n}$. 
So by $\Pi^{0}_{n}$-reflection this holds of some $ J_{\bar \theta}^{F_{\kappa}} = Q^{M}(\kappa) $ for some $\kappa < \bar \theta < \theta^{N}(\kappa)$, some $M\in H_{\kappa}^{K^{N}}$.\\ \mb
\qed

%{\em Question 1:} Are there further uses for the notions of $\Sigma_{n}$-stable measurability for $n>1$? (As this is equivalent to {\em iterability} so the question can be asked bout that too.)
%\cite{BCFHR18} 

{\em Question 1:} Under the assumptions of Lemma \ref{1.11}, does $\kappa$ inaccessible and $ u_{2}(\kappa)=\sigma(\kappa)$ imply that $\kappa$ is stably measurable?

We conjecture not, but if so, then a non-$V=K$ version of Theorem \ref{2.3} would be provable. The next question is not directly related to stable measurability but to $u_{2}$ being as large as possible. (Recall that $u_{2}(\omega_{1})$ can be $\omega_{2}$.)

{\em Question 2:} For $\kappa > \omega_{1}$ a regular cardinal, can 
$u_{2}(\kappa) = \kappa^{+}$?\\
We conjecture no.
%\newpage
\small
\bibliographystyle{plain}
\bibliography{settheory10s}

\begin{thebibliography}{10}

\bibitem{BCFHR18}
N.~Barton, A.~Caicedo, G.~Fuchs, J.D. Hamkins, and J.~Reitz.
\newblock Inner model reflection principles.
\newblock {\em ArXiv}, arXiv:1708.06669, April 2018.

\bibitem{D}
A.~J. Dodd.
\newblock {\em The {Core Model}}, volume~61 of {\em London Mathematical Society
  Lecture Notes in Mathematics}.
\newblock Cambridge University Press, Cambridge, 1982.

\bibitem{DJK}
H-D. Donder, R.~B. Jensen, and B.~Koppelberg.
\newblock Some applications of {$K$}.
\newblock In R.~Jensen and A.~Prestel, editors, {\em Set theory and Model
  theory}, volume 872 of {\em Springer Lecture Notes in Mathematics}, pages
  55--97. Springer Verlag, 1981.

\bibitem{FMW}
Q.~Feng, M.~Magidor, and W.H.Woodin.
\newblock Universally {Baire} sets of reals.
\newblock In H.~Judah, W.~Just, and W.~H. Woodin, editors, {\em Set Theory of
  the Continuum}, MSRI Publications. Springer Verlag, 1992.

\bibitem{Je72}
R.~B. Jensen.
\newblock The fine structure of the constructible hierarchy.
\newblock {\em Annals of Mathematical Logic}, 4:229--308, 1972.

\bibitem{Ka03}
A.~Kanamori.
\newblock {\em The Higher Infinite}.
\newblock Springer Monographs in Mathematics. Springer Verlag, New York, 2nd
  edition, 2003.

\bibitem{Luecke18}
P.~L\"ucke.
\newblock Partition properties for simply defined colourings.
\newblock {\em ArXive}, June 2018.

\bibitem{LSS18}
P.~L\"ucke, R-D Schindler, and P.~Schlicht.
\newblock ${\Sigma}_1 (\kappa)$-definable subsets of {$H(\kappa^+)$}.
\newblock {\em Journal for Symbolic Logic}, 82(3):1106--1131, 2017.

\bibitem{LS18}
P.~L\"ucke and P.~Schlicht.
\newblock Measurable cardinals and good {$\Sigma_1(\kappa)$}-wellorderings of
  {$H(\kappa^+)$}.
\newblock {\em Mathematical Logic Quarterly}, to appear.

\bibitem{Mi79}
W.~J. Mitchell.
\newblock Ramsey cardinals and constructibility.
\newblock {\em J. Symbolic Logic}, 44(2):260--266, 1979.

\bibitem{Sharpe-Welch2012}
I.~Sharpe and P.D. Welch.
\newblock Greatly {Erd\H{o}s} cardinals and some generalizations to the {Chang}
  and {Ramsey} properties.
\newblock {\em Annals of Pure and Applied Logic}, 162:863--902, 2011.

\bibitem{We86}
P.D. Welch.
\newblock Doing without determinacy - aspects of inner models.
\newblock In F.~Drake and J.~Truss, editors, {\em Proceedings of the Logic
  Colloquium Hull `86}, Series in Logic and its Applications, pages 333--342,
  Amsterdam, 1988. North-Holland Publishing Co.

\bibitem{W11}
P.D. Welch.
\newblock Some descriptive set theory and core models.
\newblock {\em Annals of Pure and Applied Logic}, 39:273--290, 1988.

\bibitem{W3}
P.D. Welch.
\newblock Characterising subsets of $\omega_1$.
\newblock {\em Journal of Symbolic Logic}, 59(4):1420--1432, 1994.

\bibitem{W12}
P.D. Welch.
\newblock On unfoldable cardinals, $\omega$-closed cardinals, and the
  beginnings of the inner model hierarchy.
\newblock {\em Archive for Mathematical Logic}, 43(4):443--458., 2004.

\bibitem{Z02}
M.~Zeman.
\newblock {\em Inner Models and Large Cardinals}, volume~5 of {\em Series in
  Logic and its Applications}.
\newblock de Gruyter, Berlin, New York, 2002.

\end{thebibliography}

\end{document}